\documentclass[10pt]{amsart}
\usepackage{mathrsfs}
\usepackage{amssymb}
\usepackage{amsmath}
\usepackage{amsthm}
\usepackage{color}
\usepackage{cite}
\usepackage[colorlinks=true]{hyperref}

\hypersetup{urlcolor=blue, citecolor=blue}
\makeatletter
\def\LaTeX{\leavevmode L\raise.42ex
    \hbox{\kern-.3em\size{\sf@size}{0pt}\selectfont A}\kern-.15em\TeX}
\makeatother

\newcommand{\BibTeX}{{\rm B\kern-.05em{\sc
          i\kern-.025emb}\kern-.08em\TeX}}

\newtheorem{theo}{Theorem}[section]

\newtheorem{lem}[theo]{Lemma}
\newtheorem{remark}[theo]{Remark}

\newtheorem{prop}[theo]{Proposition}

\newcommand\ee{\end{equation}}
\newcommand\bes{\begin{eqnarray}}
\newcommand\ees{\end{eqnarray}}
\newcommand\bess{\begin{eqnarray*}}
\newcommand\eess{\end{eqnarray*}}

\begin{document}
{\title[Asymptotic behavior and sharp estimates for a cooperative system]
{ Asymptotic behavior and sharp estimates for spreading fronts in a cooperative system with free boundaries}
}

\author[Q. Qin, J.J. Jiao, Z.G Wang and H. Nie]{Qian Qin$^{\dag}$, JinJing Jiao$^{\dag}$, Zhiguo Wang$^{\dag,\ast}$, Hua Nie$^{\dag,\ast}$}

\thanks{$^\dag$ School of Mathematics and Statistics, Shaanxi Normal University, Xi'an, Shaanxi {\rm710119,} China.}

\thanks{$^\ast$ The corresponding authors. E-mail addresses: niehua@snnu.edu.cn; zgwang@snnu.edu.cn.}



\makeatletter
\@namedef{subjclassname@2020}{\textup{2020} Mathematics Subject Classification}
\makeatother

\begin{abstract} \small
   This paper investigates the dynamics of a reaction-diffusion system with two free boundaries, modeling the invasion of two cooperative species, where the free boundaries represent expanding fronts. We first analyze the long-term behavior of the system, showing that it follows a spreading-vanishing dichotomy: the two species either spread across the entire region or eventually die out. In the case of spreading, we determine the asymptotic spreading speed of the fronts by using a semi-wave system and provide sharp estimates for the moving fronts. Additionally, we show that the solution to the system converges to the corresponding semi-wave solution as time tends to infinity. These results contribute to a deeper understanding of the long-term dynamics of cooperative species in reaction-diffusion systems with free boundaries.
\end{abstract}
\maketitle

\section{Introduction}\label{section 1}
\setcounter{equation}{0}
The spread of invasive species is an important issue in ecology. Numerous studies have been conducted to  better understand the nature of this spread (see, e.g., \cite{SK1997}). Notably, Du and Lin \cite{YD2010} highlighted that Skellam \cite{S1951} was among the first to observe that the spread  of muskrats in Europe in the early 1900s followed a linear pattern, with the spreading radius increasing linearly over time.  This pattern has also been observed in other species \cite{SK1997}. The modeling of biological invasion has garnered considerable attention from both ecologists and mathematicians, resulting in significant insights into species spreading through various mathematical models (see, e.g., \cite{ANS2015-1}, \cite{AW1978},\cite{Fisher1937},\cite{KPP1937} and references therein).

For single-species models, one of the most successful mathematical frameworks is the diffusive logistic equation in $\mathbb{R}$, introduced by Fisher \cite{Fisher1937} and Kolmogorov,  Petrovsky,  and   Piskunov (KPP) \cite{KPP1937} to study the spread of an advantageous gene in a population.
Fisher proposed that there exists a spreading speed for the new gene, corresponding to the minimal spreading speed of an associated traveling wave problem. This proposition was convincingly proved by Aronson and Weinberger  \cite{AW1978}. Although the Fisher-KPP model effectively determines the spreading speed, it does not accurately locate the spreading front \cite{Du2022BMS}. To address this limitation, Du and Lin \cite{YD2010} introduced a free boundary version of the Fisher-KPP equation, where the species' habitat is represented by a changing interval $(g(t),h(t))$ that adheres to a well-known Stefan condition. They demonstrated that the dynamics of the solution are characterized by a spreading-vanishing dichotomy. Subsequently, Du, Matsuzawa, and Zhou \cite{DMZ2014} further refined the understanding of the sharp long-term behavior of the solution. For other single-species models with free boundaries, sharp convergence results have been achieved in several works; see \cite{DMZ2015},  \cite{KMY2020} and references therein.

In natural environments, multiple species coexist and interact, significantly influencing each other's dynamics. To accurately model the dynamics, it is essential to consider interspecies interactions.  Cooperation, where species mutually benefit from their interactions, is a prevalent form of interaction in ecological systems. To model the spread of gonorrhea, Lajmanovich and Yorke \cite{LajmanovichYorke1976} analyzed  an ODE system involving $n$ cooperative species. Subsequently,
Molina-Meyer   \cite{MM1996} investigated the following reaction-diffusion  system for two cooperative species:
\begin{equation}\label{Cauchy Problem}
\begin{cases}
u_{t}-d_1\Delta u=-au+bv-F(u), & x\in\Omega, \ t>0,\\
v_{t}-d_2\Delta v=cu-dv-G(v), & x\in\Omega, \ t>0,\\
u(x,t)=v(x,t)=0,& x\in\partial\Omega,\ t>0,\\
u(x,0)=u_0(x),\ v(x,0)=v_0(x),& x\in\Omega,
\end{cases}
\end{equation}
where $\Omega$ is a bounded domain of $\mathbb{R}^N~(N\ge1)$ with a smooth boundary $\partial\Omega$, and $d_1,d_2,a,b,c$, and $d$ are positive constants.
These kinds of systems are prevalent in various fields, including chemistry, ecology, and epidemiology. They can describe the evolution of two cooperative species with densities $u(x,t)$ and $v(x,t)$. Here, $d_1$ and $d_2$ are the diffusion rates of the species, $b$ and $c$ measure the mutual cooperative effects, and $a$ and $d$ represent the mortality rates. The functions $F(u)$ and $G(v)$ characterize the population loss of the species and are assumed to  satisfy the following conditions:
\begin{itemize}
  \item[$\mathrm{(A1)}$]~$F,G\in C^2([0,+\infty))$ are non-decreasing and $F(0)=G(0)=F'(0)=G'(0)=0$;
  \item[$\mathrm{(A2)}$]~$F(\zeta)/\zeta$ and $G(\zeta)/\zeta$ are strictly increasing with respect to $\zeta\in(0,+\infty)$;
  \item[$\mathrm{(A3)}$]~$\lim\limits_{\zeta\to+\infty}F(\zeta)/\zeta=F'(+\infty)=+\infty$ and $\lim\limits_{\zeta\to +\infty}G(\zeta)/\zeta=G'(+\infty)=+\infty$.
\end{itemize}

In   \cite{MM1996}, Molina-Meyer focused on the positive steady state of system \eqref{Cauchy Problem}, examining its existence, uniqueness, and stability. Numerous studies have explored the dynamics of system \eqref{Cauchy Problem} with  specific choices of $F(u)$ and $G(v)$. For instance, \'{A}lvarez-Caudevilla and L\'{o}pez-G\'{o}mez  \cite{AL2010} analyzed the impacts of cooperative parameters $b$ and $c$ on the dynamics when $G(v)\equiv0$. \'{A}lvarez-Caudevilla, Du, and Peng  \cite{ADP2014} studied the long-term behavior of the solution when $F(u)=u^p$ and $G(v)=v^q$ $(p,q>1)$. For additional research works related to system \eqref{Cauchy Problem}, see  \cite{CM1981},  \cite{ZW2004} and their references.

To better locate the spreading fronts of the two cooperative species,  it is practical to incorporate dynamic changes in their living regions over time, as suggested by Du and Lin \cite{YD2010}. By taking  $N=1$ and  including time-dependent variations in the living regions, we derive the following cooperative reaction-diffusion system with two free boundaries:
\begin{equation}\label{1.2}
\begin{cases}
u_{t}-d_1 u_{xx}=-au+bv-F(u), & g(t)<x<h(t), \ t>0,\\
v_{t}-d_2 v_{xx}=cu-dv-G(v), & g(t)<x<h(t), \ t>0,\\
u(h(t),t)=v(h(t),t)=u(g(t),t)=v(g(t),t)=0,& t>0,\\
 h'(t)=-\mu[v_{x}(h(t),t)+\rho u_{x}(h(t),t)],& t>0,\\
g'(t)=-\mu[v_{x}(g(t),t)+\rho u_{x}(g(t),t)],& t>0,\\
g(0)=-h_0,~h(0)=h_0,~u(x,0)=u_0(x),\ v(x,0)=v_0(x),& -h_0\leq x\leq h_0.
\end{cases}
\end{equation}
Here, $u(x,t)$ and $v(x,t)$ represent the population densities of the two cooperative species at location $x$ and time $t$.
The parameters  $d_1,d_2,a,b,c,d$  and functions $F(u)$, $G(v)$ retain their previous meanings from system  \eqref{Cauchy Problem}, with $F$ and $G$ satisfying  the conditions $\mathrm{(A1)}$-$\mathrm{(A3)}$.
The  positive constants $\mu$ and $\rho$ quantify the species' capacity to invade  new territories.
The functions $x=h(t)$ and  $x=g(t)$ denote the moving boundaries  to be determined, which adhere to the well-known Stefan boundary conditions (the fourth and fifth equations in \eqref{1.2}).  This Stefan condition indicates that the expanding rate of the region $(g(t),h(t))$ is proportional to a linear combination of the spatial gradients of $u$ and $v$ at the boundaries.
To ensure a classical solution for system  \eqref{1.2}, we assume that the initial functions $u_0$ and $v_0$ are defined on the interval $[-h_0,h_0]$ and satisfy
\begin{equation}
\begin{cases}
u_0\in C^{2}([-h_0,h_0]),~u_0(\pm h_0)=0,~\textrm{and}~u_0(x)>0\text{ for }x\in (-h_0,h_0),\\
v_0\in C^{2}([-h_0,h_0]),~v_0(\pm h_0)=0,~\textrm{and}~v_0(x)>0\text{ for }x\in (-h_0,h_0).
\end{cases}\label{1.3}
\end{equation}

The primary objective of this paper is to investigate the dynamics of system \eqref{1.2} with a  particular focus on the long-term behavior of its solutions.
Given the cooperative nature of system \eqref{1.2}, the comparison principle remains applicable. Our approach involves constructing upper and lower solutions to explore its dynamical properties, though the nonlinearities in $F(u)$ and $G(v)$, along with the interaction of two species, complicate this process (refer to Lemma \ref{lem2.6}). Similar to the free boundary version of the Fisher-KPP equation discussed in \cite{YD2010}, system \eqref{1.2} exhibits a spreading-vanishing dichotomy (see Theorem \ref{theo1.1}). Specifically, this result implies global convergence within the region $(g(t),h(t))$ when spreading does not occur.

The second objective of this paper is to analyze the asymptotic spreading speeds of the free boundaries $h(t)$ and $g(t)$ as $t\to+\infty$ in scenarios where spreading occurs. Additionally, we aim to gain a deeper understanding of the global  behavior of the population densities $u(x,t)$ and $v(x,t)$.  To achieve these goals, our approach proceeds in  two steps: (1) Establish asymptotic spreading speeds. Drawing inspiration from  Wang, Nie, and Du \cite{WND2019}, we  first derive a preliminary result regarding the asymptotic spreading speeds, given by $\lim\limits_{t\to+\infty}h(t)/t$ and $-\lim\limits_{t\to+\infty}g(t)/t$ (see Theorem \ref{theo1.3}). This involves employing upper and lower solution methods alongside a perturbation argument. However, since the Stefan boundary conditions (the fourth and fifth equations in \eqref{1.2}) differ significantly  from those in \cite{WND2019}, the derivation requires substantial modifications (see Subsection \ref{subsection 4.1}).
(2) Derive   sharp estimates and global profiles. In the second step, we seek to obtain sharper estimates for the free boundaries  and   global profiles of the population densities (see Theorem \ref{theo1.4}). This analysis is inspired by works \cite{DMZ2014},  \cite{WND},  \cite{WQW2022}, where \cite{DMZ2014} addresses single-species equations, \cite{WND} examines a West Nile virus system, and \cite{WQW2022} investigates a Lotka-Volterra  competition system. Since system \eqref{1.2} involves two species and the Stefan boundary conditions differ  from those in \cite{WND},  \cite{WQW2022}, we encounter two main difficulties: (i) The calculations to derive bounds for $h(t)+s_{\mu,\rho} t$ and $g(t)-s_{\mu,\rho} t$, where $s_{\mu,\rho}$ is a crucial propagation speed for an associated semi-wave problem, are more intricate (see Subsection \ref{subsection 4.2}). (ii) Demonstrating  that the solution to  system \eqref{1.2} converges to the solution of an associated semi-wave system poses a nontrivial challenge (refer to Section \ref{section 5}).

Before ending this section, we would like to highlight additional references that form part of the background of this research. Inspired by the work of Du and Lin \cite{YD2010}, substantial efforts have been made to develop analytical tools for more general single-species models with free boundaries. Notable studies include \cite{DGP2013}, which addresses problems in periodic environments, and \cite{2017DingPengWeiJDE}, which focuses on advection-related issues. Extensive research on multiple species has produced interesting results, as documented in \cite{2015GuoWuNonlinearity}, \cite{2023LiWuANS}, \cite{2007LiangZhaoCPAM}, \cite{ANS2017-1},  \cite{2014WangZhaoJDDE},  \cite{2017WangZhaoJDDE},  \cite{WNND2019}.  Specifically,
\cite{ANS2017-1}  explores cooperative systems, \cite{2015GuoWuNonlinearity}, \cite{2014WangZhaoJDDE} investigate competition systems,   \cite{2023LiWuANS},
\cite{2017WangZhaoJDDE} examine prey-predator systems, and
\cite{2007LiangZhaoCPAM},  \cite{WNND2019} focus on West Nile virus-related problems.

The rest of this paper is organized as follows. In Section \ref{section 2}, we present the main results of the paper. In Section \ref{section 3}, we provide preliminary results for  system \eqref{1.2}, including the well-posedness and comparison principles, followed by the derivation of the spreading-vanishing dichotomy (Theorem \ref{theo1.1}). Section  \ref{section 4} is devoted to analyzing the asymptotic spreading speeds of the free boundaries (Theorem \ref{theo1.3}). In Section  \ref{section 5}, we derive  sharp estimates for the free boundaries and global profiles of the population densities (Theorem \ref{theo1.4}), based on the estimates from Section  \ref{section 4} and the construction of suitable upper and lower solutions. Finally,  a brief discussion concludes the paper in Section \ref{section 6}.

\section{Main results}\label{section 2}
\setcounter{equation}{0}
The main findings of this paper are presented in this section. Before we proceed, it is important to mention the following two points:
\begin{itemize}
\item [$\bullet$] If $\mathrm{(A1)}$-$\mathrm{(A3)}$ hold and $(u_0,v_0)$ satisfies \eqref{1.3}, then system \eqref{1.2}  admits a unique  solution $(u,v,g,h)=(u(x,t),v(x,t),g(t),h(t))$ (see Theorem \ref{lem2.0} in Section \ref{section 3}).
\item [$\bullet$] If  $\mathrm{(A1)}$-$\mathrm{(A3)}$ hold and $bc-ad>0$,  then the  ODE system  corresponding to \eqref{1.2}  has a unique positive equilibrium $(u^*,v^*)$ that is globally asymptotically stable (see, e.g., \cite{MM1996}).
\end{itemize}
Throughout the paper, we make the following assumptions:
\begin{itemize}
	\item [(H)] \quad $\mathrm{(A1)}$-$\mathrm{(A3)}$ hold, $(u_0,v_0)$ satisfies \eqref{1.3},  and  $bc-ad>0$.
\end{itemize}

The following results present a  spreading-vanishing dichotomy for system \eqref{1.2}.
\begin{theo}\label{theo1.1}
{\rm(Dichotomy)}
Assume that $\mathrm{(H)}$ holds and let  $(u,v,g,h)$ be the unique solution of system \eqref{1.2}. Then the following two limits are well-defined:
$$h_{\infty}:=\lim_{t\to+\infty}h(t)\text{ and }g_{\infty}:=\lim_{t\to+\infty}g(t).$$
 Moreover, one of the following cases occurs:
\begin{itemize}
  \item[(i)]Vanishing:
  $h_{\infty}-g_{\infty}<+\infty$ and
  $\lim\limits_{t\to+\infty}(u(\cdot,t),v(\cdot,t))=(0,0)$ uniformly in $[g(t),h(t)]$;
  \item[(ii)]Spreading:
  $h_{\infty}=-g_{\infty}=+\infty$ and $\lim\limits_{t\to+\infty}(u(\cdot,t),v(\cdot,t))
  =(u^*,v^*)$ locally uniformly in $\mathbb{R}$.
\end{itemize}
\end{theo}

From a mathematical perspective, Theorem \ref{theo1.1}  describes the long-term  behavior of the solution to system \eqref{1.2}, where convergence is global in the vanishing case.
Biologically, this means that the two cooperative species will either  eventually die out (vanishing) or spread throughout the entire region (spreading). In the spreading scenario, the two species will achieve densities  $u^*$ and $v^*$ respectively in any bounded domain of $\mathbb{R}$. To further explore the spreading dynamics, two natural questions arise:

\emph{Question 1}: Can the asymptotic spreading speed of the fronts $h(t)$ and $g(t)$ be determined?

\emph{Question 2}:  Is it possible to determine the global behavior of the densities $u(x,t)$ and $v(x,t)$?

To answer the two questions above, we need to investigate the following  semi-wave system:
\begin{equation}\label{1.5}
\begin{cases}
d_1 \phi''-s\phi'-a\phi+b\psi-F(\phi)=0, &0<\xi<+\infty,\\
d_2 \psi''-s\psi'+c\phi-d\psi-G(\psi)=0, &0<\xi<+\infty,\\
(\phi(0),\psi(0))=(0,0),~(\phi(+\infty),\psi(+\infty))=(u^*,v^*),
\end{cases}
\end{equation}
where $bc-ad>0$, and $(u^*,v^*)$ is the unique positive equilibrium of the  ODE system corresponding to \eqref{1.2}.

\begin{prop}\label{theo1.2}
{\rm(Semi-wave)} Assume that $\mathrm{(A1)}$-$\mathrm{(A3)}$ hold and  $bc-ad>0$. Then there exists a constant $s^{\ast}>0$ such that for every   $s\in[0,s^{\ast})$, system \eqref{1.5} has a unique strictly increasing solution $(\phi_{s},\psi_{s})\in C^{2}([0,+\infty))\times C^{2}([0,+\infty))$, and system \eqref{1.5} has no such solution when  $s\ge s^{\ast}$. Moreover, for any $\mu>0$ and $\rho>0$, there exists a unique $s_{\mu,\rho}\in(0,s^{\ast})$ such that
\begin{equation}\label{s mu}
\mu[\psi'_{s_{\mu,\rho}}(0)+\rho\phi'_{s_{\mu,\rho}}(0)]=s_{\mu,\rho}.
\end{equation}
\end{prop}

With the help of Proposition $\ref{theo1.2}$, we obtain the asymptotic spreading speeds of the fronts $h(t)$ and $g(t)$ when spreading occurs. The results are presented as follows.

\begin{theo}\label{theo1.3}
{\rm(Spreading speed)} Assume that $\mathrm{(H)}$ holds and let $(u,v,g,h)$ be the unique solution of system \eqref{1.2}. If spreading occurs, then
\begin{equation}\label{1.6}
\lim_{t \to+\infty}\frac{h(t)}{t}=-\lim_{t \to+\infty}\frac{g(t)}{t}=s_{\mu,\rho},
\end{equation}
where $s_{\mu,\rho}\in(0,s^*)$ is uniquely determined by \eqref{s mu}.
\end{theo}

 Theorem \ref{theo1.3} suggests  that, in the spreading scenario, the moving fronts $h(t)$ and $g(t)$ share the same asymptotic spreading speed $s_{\mu,\rho}$, which reflects the propagation speed of the semi-wave system \eqref{1.5} influenced by species' invasion ability (see \eqref{s mu}). Moreover, it follows from \eqref{1.6} that
$$h(t)=s_{\mu,\rho} t+o(t)\text{ and }g(t)=-s_{\mu,\rho} t+o(t)\text{ as }t\to+\infty.$$

By analyzing the semi-wave system \eqref{1.5}, we can further determine  $o(t)$ and obtain a global description of the population densities $u(x,t)$ and $v(x,t)$ during the spreading process.
\begin{theo}\label{theo1.4}
{\rm(Sharp profile)} Assume that $\mathrm{(H)}$ holds, and let $(u,v,g,h)$ and $(\phi_{s_{\mu,\rho}},\psi_{s_{\mu,\rho}})$ be the solutions of systems \eqref{1.2} and \eqref{1.5}, respectively. If spreading occurs, then there exist constants $g^{\ast}$ and $h^{\ast}$ such that
\begin{align}
\lim_{t \to +\infty}(g(t)+s_{\mu,\rho}t-g^{\ast})=0,\ \lim_{t \to +\infty}g'(t)=-s_{\mu,\rho},\label{1.7}\\
\lim_{t \to +\infty}(h(t)-s_{\mu,\rho}t-h^{\ast})=0,\ \lim_{t \to +\infty}h'(t)=s_{\mu,\rho}.\label{1.8}
\end{align}
Moreover,
\begin{align}
\lim_{t \to +\infty}\|(u(\cdot,t),v(\cdot,t))-
(\phi_{s_{\mu,\rho}}(\cdot-g(t)),\psi_{s_{\mu,\rho}}(\cdot-g(t)))\|_{[L^{\infty}([g(t),0])]^2}=0,\label{1.9}\\
\lim_{t \to +\infty}\|(u(\cdot,t),v(\cdot,t))-
(\phi_{s_{\mu,\rho}}(h(t)-\cdot),\psi_{s_{\mu,\rho}}(h(t)-\cdot))\|_{[L^{\infty}([0,h(t)])]^2}=0.\label{1.10}
\end{align}
\end{theo}

Theorem \ref{theo1.4} shows that the long-term behavior of the solution  $(u,v,g,h)$ to system  \eqref{1.2} can be comprehensively understood via semi-wave system \eqref{1.5}.  To be more specific, \eqref{1.7} and \eqref{1.8} provide a more precise characterization of the moving fronts $h(t)$ and $g(t)$ compared to \eqref{1.6}, revealing that  these spreading fronts asymptotically converge to linear functions with constant correction terms. Additionally, it follows from  \eqref{1.9} and \eqref{1.10} that the global profiles of the population densities $u(x,t)$ and $v(x,t)$  within the domain $(g(t),h(t))$ can be described by the semi-wave solution $(\phi_{s_{\mu,\rho}},\psi_{s_{\mu,\rho}})$ of system \eqref{1.5}. While there is substantial  literature on spreading-vanishing dichotomies and spreading speeds in various multi-species evolutionary systems (see, e.g., \cite{WD2021}, \cite{WNND2019}) as well as sharp profiles for single-species models (see, e.g., \cite{DMZ2015}, \cite{KMY2020},  \cite{ANS2023-1}), studies specifically addressing the sharp asymptotic spreading speed of fronts and global population profiles remain relatively limited.
Thus,  Theorem \ref{theo1.4} represents a significant contribution, offering precise dynamical insights in this context.

\section{Spreading-vanishing dichotomy}\label{section 3}
\setcounter{equation}{0}
In this section, we first present  the well-posedness  and  comparison principles for system \eqref{1.2}. Subsequently, we proceed with the proof of Theorem \ref{theo1.1} by constructing appropriate upper and lower solutions.
For clarity, we introduce the following notations, which will be used throughout this paper. For any vectors $\mathbf{p}=(p_1,\cdots,p_m)$ and $\mathbf{q}=(q_1,\cdots,q_m)$ in $\mathbb{R}^{m}$,
$\mathbf{p}\preceq(\succeq)\mathbf{q}$ (resp. $\mathbf{p}\prec(\succ)\mathbf{q})$ denotes that  $p_{i}\leq(\geq)~q_{i}$ (resp. $p_{i}<(>)~q_{i})$ for $i=1,\cdots,m$.
The transpose of a matrix $B$ is denoted by $B^{\intercal}$.

The following theorem establishes the global existence and uniqueness of a positive solution to system \eqref{1.2}. For a detailed proof, please refer to \cite[Theorem 1.1]{WD2021}.
\begin{theo}\label{lem2.0}
{\rm(Existence and uniqueness)} Assume that  $\mathrm{(A1)}$-$\mathrm{(A3)}$ hold and $(u_0,v_0)$ satisfies \eqref{1.3}.  For any constant $\alpha\in(0,1)$, system \eqref{1.2} admits a unique solution
$$(u,v,g,h)\in[C^{2+\alpha,1+\frac{\alpha}{2}}(D)]^{2}\times [C^{1+\frac{1+\alpha}{2}}((0,+\infty))]^2,$$
where $D:=\{(x,t)\in\mathbb{R}^{2}:x\in[g(t),h(t)],~t>0\}$. Furthermore, there exists a constant $C_{0}>0$ such that
\begin{align*}
&0<u(x,t),~v(x,t)\le C_0\text{ for }x \in (g(t),h(t)),~t>0,\\
&0<-g'(t),~h'(t) \le C_0\text{ for }t>0.
\end{align*}
\end{theo}

 We next present two comparison principles for system \eqref{1.2}. They are straightforward variations of  \cite[Lemmas 2.1 and 2.2]{WND} and can be proved using arguments similar to those in \cite[Lemma 2.6]{DL2014}.

\begin{lem}\label{lem2.1}
 Let $(u,v,g,h)$ be the solution of  \eqref{1.2}. Assume that $T\in(0,+\infty)$, $\overline{g},\overline{h} \in C^{1}([0,T])$, $g(t)\leq \overline{g}(t)<\overline{h}(t)$ in
$[0,T]$,
$\overline{u}, \overline{v} \in C(\overline{D_{T}})\cap C^{2,1}(D_{T})$ with  $D_{T}=\{(x,t) \in \mathbb{R}^{2}: x \in (\overline{g}(t),\overline{h}(t)), t \in (0,T]\}$, and
\begin{equation}\label{2.1}
\begin{cases}
    \overline{u}_{t}-d_1 \overline{u}_{xx}\geq-a\overline{u}+b\overline{v}-F(\overline{u}), & \overline{g}(t)<x<\overline{h}(t), \ 0<t<T,\\
    \overline{v}_{t}-d_2 \overline{v}_{xx}\geq c\overline{u}-d\overline{v}-G(\overline{v}), & \overline{g}(t)<x<\overline{h}(t), \ 0<t<T,\\
    \overline{u}(x,t)\geq u(x,t),~\overline{v}(x,t)\geq v(x,t), & x=\overline{g}(t),~0<t<T,\\
    \overline{u}(x,t)=\overline{v}(x,t)=0, & x=\overline{h}(t),~ 0<t<T,\\
    \overline{h}'(t)\geq -\mu[\overline{v}_{x}(\overline{h}(t),t)+\rho\overline{u}_{x}(\overline{h}(t),t)],& 0<t<T,\\
    \overline{u}(x,0)\geq u_0(x),~\overline{v}(x,0)\geq v_0(x),& \overline{g}(0)\leq x\leq h_0.
\end{cases}
\end{equation}
 Then
$$h(t)\leq\overline{h}(t),~u(x,t)\leq\overline{u}(x,t),~v(x,t)\leq\overline{v}(x,t)\text{ for }
\overline{g}(t)\le x\le h(t),~0<t\le T.$$
\end{lem}

\begin{lem}\label{lem2.1.2}
 Let $(u,v,g,h)$ be the solution of system \eqref{1.2}. Assume that $T\in(0,+\infty)$,
$\overline{g},\overline{h} \in C^{1}([0,T])$, $\overline{g}(t)<\overline{h}(t)$ in
$[0,T]$,
$\overline{u}, \overline{v} \in C(\overline{D_{T}})\cap C^{2,1}(D_{T})$ with $D_{T}=\{(x,t) \in \mathbb{R}^{2}: x \in (\overline{g}(t),\overline{h}(t)),~t \in (0,T]\}$, and
\begin{equation}\label{2.1.2}
\begin{cases}
    \overline{u}_{t}-d_1 \overline{u}_{xx}\geq-a\overline{u}+b\overline{v}-F(\overline{u}), & \overline{g}(t)<x<\overline{h}(t), \ 0<t<T,\\
    \overline{v}_{t}-d_2 \overline{v}_{xx}\geq c\overline{u}-d\overline{v}-G(\overline{v}), & \overline{g}(t)<x<\overline{h}(t), \ 0<t<T,\\
    \overline{u}(x,t)=\overline{v}(x,t)=0, & x=\overline{g}(t)\text{ or }x=\overline{h}(t),~0<t<T,\\
    \overline{g}'(t)\leq -\mu[\overline{v}_{x}(\overline{g}(t),t)+\rho\overline{u}_{x}(\overline{g}(t),t)],& 0<t<T,\\
    \overline{h}'(t)\geq -\mu[\overline{v}_{x}(\overline{h}(t),t)+\rho\overline{u}_{x}(\overline{h}(t),t)],& 0<t<T,\\
    \overline{u}(x,0)\geq u_0(x),~\overline{v}(x,0)\geq v_0(x),&-h_0\leq x\leq h_0.
\end{cases}
\end{equation}
Then
$$g(t)\ge\overline{g}(t),~h(t)\leq\overline{h}(t),~u(x,t)\leq\overline{u}(x,t),~v(x,t)\leq\overline{v}(x,t)
\text{ for }g(t)\le x\le h(t),~0<t\le T.$$
\end{lem}

\begin{remark}\label{remark2.2}
We present the following comments regarding Lemmas \ref{lem2.1} and \ref{lem2.1.2}.
\begin{itemize}
		\item[$\mathrm{(i)}$] If the reversed inequalities in \eqref{2.1} hold, and $(\overline{u},\overline{v},\overline{g},\overline{h})$ is rewritten as $(\underline{u},\underline{v},\underline{g},\underline{h})$, then
$$h(t)\geq\underline{h}(t),~u(x,t)\geq\underline{u}(x,t),~v(x,t)\geq\underline{v}(x,t)\text{ for } g(t)\le x\le\underline{h}(t),~0<t\le T.$$
\item[$\mathrm{(ii)}$] Similarly, if the reversed inequalities in \eqref{2.1.2} hold, and $(\overline{u},\overline{v},\overline{g},\overline{h})$ is rewritten as $(\underline{u},\underline{v},\underline{g},\underline{h})$, then
$$g(t)\le\underline{g}(t),~h(t)\geq\underline{h}(t),~u(x,t)\geq\underline{u}(x,t),~v(x,t)\geq\underline{v}(x,t)\text{ for } \underline{g}(t)\le x\le \underline{h}(t),~0<t\le T.$$
\item[$\mathrm{(iii)}$] The functions $(\overline{u},\overline{v},\overline{g},\overline{h})$ and $(\underline{u},\underline{v},\underline{g},\underline{h})$ are usually called an upper solution and a lower solution of system \eqref{1.2}, respectively.
\end{itemize}
\end{remark}
The strict monotonicity of $h(t)$ and $g(t)$ in  Theorem \ref{lem2.0} implies that $h_{\infty}\in(h_0,+\infty]$ and $g_{\infty}\in[-\infty,-h_0)$. If $h_{\infty}<+\infty$ or $g_{\infty}>-\infty$, then both $h_{\infty}$ and $g_{\infty}$ are finite (see, e.g., \cite[Lemma 5.2]{YD2010}). Therefore, one of the following alternatives must hold: (i)  $h_{\infty}-g_{\infty}<+\infty$, or (ii)  $h_{\infty}=-g_{\infty}=+\infty$. With this in mind, we proceed to prove the spreading-vanishing dichotomy (i.e., Theorem  \ref{theo1.1}). For clarity, the proof is divided into the following two lemmas.
\begin{lem}\label{lem2.5}
Assume that $\mathrm{(H)}$ holds and let $(u,v,g,h)$ be the unique solution of system \eqref{1.2}.
If $h_{\infty}-g_{\infty}<+\infty$, then $\lim\limits_{t\to+\infty}\|u(\cdot,t)\|_{C([g(t),h(t)])}=\lim\limits_{t\to+\infty}\|v(\cdot,t)\|_{C([g(t),h(t)])}=0$.
\end{lem}
\begin{proof}
Since the proofs of the two limits are similar, we only prove that $\lim\limits_{t\to+\infty}\|v(\cdot,t)\|_{C([g(t),h(t)])}=0$.
It suffices to prove that $\limsup\limits_{t\to+\infty}\|v(\cdot,t)\|_{C([g(t),h(t)])}=0$.
Otherwise, since $v(x,t)$ is uniformly bounded in $[g(t),h(t)]\times(0,+\infty)$ (see Theorem \ref{lem2.0}), there  exists $\delta_*>0$ such that $\limsup\limits_{t\to+\infty}\|v(\cdot,t)\|_{C([g(t),h(t)])}\geq\delta_*$. Then one can find a sequence $(x_{k},t_{k})\in[g(t_{k}),h(t_{k})]\times(0,+\infty)$ with $\lim\limits_{k\to+\infty}t_k=+\infty$ such that $v(x_{k},t_{k})\geq\delta_*/2$ for all $k\in\mathbb{N}$. By Theorem \ref{lem2.0}, we have $h(t_{k})-g(t_{k})<h_\infty-g_\infty<+\infty$ for any $k\in\mathbb{N}$, which implies that $\{x_{k}\}$ is bounded. Thus, there exists a subsequence of $\{x_{k}\}$ (still denoted it  by $\{x_{k}\}$) such that $\lim\limits_{k\to+\infty}x_k=x_0\in[g_{\infty},h_{\infty}]$.
By similar arguments as in the proofs of \cite[Theorems 1.1 and 2.1]{WD2021}, one can find a constant $M>0$ such that
$$\Vert h'\Vert_{C^{\frac{\alpha}{2}}([\tau,\tau+1])}+
\Vert g'\Vert_{C^{\frac{\alpha}{2}}([\tau,\tau+1])}\le M$$
 for any $\tau\ge1$.
Therefore, since $\lim\limits_{t\to+\infty}h(t)=h_{\infty}$ is finite, we can conclude that $\lim\limits_{t\to+\infty}h'(t)=0$. It then follows from the fourth equation in \eqref{1.2} that
\begin{equation}\label{limit01}
\lim\limits_{k\to+\infty}[v_{x}(h(t_{k}),t_{k})+\rho u_{x}(h(t_{k}),t_{k})]=0.
\end{equation}

Clearly, there exists $t_0>0$  such that $|h(t)-h(t_0)|<h(t_0)/12$ and $|g(t)-g(t_0)|<-g(t_0)/12$ for $t\ge t_0$.
Choose $\beta, \gamma\in C^3(\mathbb{R})$ such that
$$\beta(y)=1\text{ if }|y-h(t_0)|<\frac{h(t_0)}{4},\  \beta(y)=0\text{ if }|y-h(t_0)|>\frac{h(t_0)}{2},\ |\beta'(y)|<\frac{5}{h(t_0)}\text{ for all }y,$$
and
$$\gamma(y)=-1\text{ if }|y-g(t_0)|<-\frac{g(t_0)}{4},\  \gamma(y)=0\text{ if }|y-g(t_0)|>-\frac{g(t_0)}{2},\ |\gamma'(y)|<-\frac{5}{g(t_0)}\text{ for all }y.$$

We next straighten the free boundaries by the method developed in \cite{ChenFriedman2000}, \cite{Wang2019DCDSB}.
Consider  the transformation $(y,t)\mapsto (x,t)$ given by
$$x=\Theta(y,t):=y+\beta(y)[h(t)-h(t_0)]-\gamma(y)[g(t)-g(t_0)]\text{ for }y\in\mathbb{R},$$
which, for every fixed $t\ge t_0$, is a diffeomorphism from $\mathbb{R}$ to $\mathbb{R}$.  Moreover, it maps $x=h(t)$ and $x=g(t)$ to  the lines $y=h(t_0)$ and $y=g(t_0)$, respectively.

Direct calculations yield that
\begin{equation}\label{add*}
\begin{split}
   \frac{\partial y}{\partial x} &=\frac{1}{1+\beta'(y)[h(t)-h(t_0)]-\gamma'(y)[g(t)-g(t_0)]}=:\sqrt{A(y,h(t),g(t))}, \\
    \frac{\partial^2 y}{\partial x^2} &=-\frac{\beta''(y)[h(t)-h(t_0)]-\gamma''(y)[g(t)-g(t_0)]}{[1+\beta'(y)(h(t)-h(t_0))-\gamma'(y)(g(t)-g(t_0))]^3}=:B(y,h(t),g(t)),\\
   \frac{\partial y}{\partial t}&=\frac{-h'(t)\beta(y)+g'(t)\gamma(y)}{1+\beta'(y)[h(t)-h(t_0)]-\gamma'(y)[g(t)-g(t_0)]}\\
   &=:-h'(t)C(y,h(t),g(t))+g'(t)D(y,h(t),g(t)).
\end{split}
\end{equation}
Denote $\hat{U}(y,t)=u(\Theta(y,t),t)$ and $\hat{V}(y,t)=v(\Theta(y,t),t)$ for $y\in[g(t_0),h(t_0)]$ and $t\ge t_0$. Define
$$u_k(y,t):=\hat{U}(y,t+t_k),\ v_k(y,t):=\hat{V}(y,t+t_k)\text{ for }y\in[g(t_0),h(t_0)],\ t\ge t_0-t_k.$$
From  the estimates for $(u,v,g,h)$ in Theorem \ref{lem2.0}, it follows that $(u_k,v_k)$ has a subsequence $(u_{k_i},v_{k_i})$ converging to some $(\hat{u},\hat{v})$ as $i\to+\infty$ in $\big[C_{\mathrm{loc}}^{2,1}([g(t_0),h(t_0)]\times\mathbb{R})\big]^{2}$, and $(\hat{u},\hat{v})$ satisfies
\begin{equation}\label{add equation 01}
\begin{cases}
\hat{u}_{t}=d_1A_{\infty}(y)\hat{u}_{yy}+d_1B_{\infty}(y)\hat{u}_{y}-a\hat{u}+b\hat{v}-F(\hat{u}), & g(t_0)<y<h(t_0), \ t\in\mathbb{R},\\
\hat{v}_{t}=d_2A_{\infty}(y)\hat{v}_{yy}+d_2B_{\infty}(y)\hat{v}_{y}+c\hat{u}-d\hat{v}-G(\hat{v}), & g(t_0)<y<h(t_0), \ t\in\mathbb{R},\\
\hat{u}(y,t)\ge0,~\hat{v}(y,t)\ge0, &  g(t_0)<y<h(t_0), \ t\in\mathbb{R},\\
\hat{u}(y,t)=\hat{v}(y,t)=0, &y=g(t_0)\text{ or }y=h(t_0),~t\in\mathbb{R},
\end{cases}
\end{equation}
where $A_{\infty}(y):=A(y,h_{\infty},g_{\infty})$ and $B_{\infty}(y):=B(y,h_{\infty},g_{\infty})$, with $A$ and $B$ defined in \eqref{add*}.

Let $y_0\in[g(t_0),h(t_0)]$ satisfy $x_0=\Theta(y_0,+\infty)$. Then $\hat{v}(y_0,0)\ge\delta_*/2$. Thus, it follows from the strong maximum principle   that $\hat{v}>0$ in $(g_{\infty},h_{\infty})\times\mathbb{R}$.
This, together with $\mathrm{(A1)}$ and \eqref{add equation 01}, implies that $\hat{u}\not\equiv0$ and hence  $\hat{u}>0$  in $(g_{\infty},h_{\infty})\times\mathbb{R}$. Applying the Hopf boundary lemma, we obtain
$$\hat{u}_y(h(t_0),t)<0\text{ and }\hat{v}_y(h(t_0),t)<0\text{ for }t\in\mathbb{R},$$
which implies that
$$\lim\limits_{k\to+\infty}[v_{x}(h(t_{k}),t_{k})+\rho u_{x}(h(t_{k}),t_{k})]=\sqrt{A_{\infty}(h(t_0))}\big[\hat{v}_{y}(h(t_0),0)+\rho\hat{u}_{y}(h(t_0),0)\big]<0,$$
contradicting \eqref{limit01}. Hence we have $\lim\limits_{t\to+\infty}\|v(\cdot,t)\|_{C([g(t),h(t)])}=0$, which  completes the proof.
\end{proof}

\begin{lem}\label{lem2.6}
Assume that $\mathrm{(H)}$ holds and let $(u,v,g,h)$ be the unique solution of system \eqref{1.2}.
If $h_{\infty}=-g_{\infty}=+\infty$, then $\lim\limits_{t\to+\infty}(u(\cdot,t),v(\cdot,t))=(u^*,v^*)$ locally uniformly in  $\mathbb{R}$.
\end{lem}
\begin{proof}
Let $(U(t),V(t))$ be the solution of following ODE system:
\begin{equation}\label{1.11}
\begin{cases}
U_t=-aU+bV-F(U), &t>0,\\
V_t=cU-dV-G(V), &t>0,\\
U(0)=\|u_{0}\|_{C([-h_0,h_0])},\ V(0)=\|v_{0}\|_{C([-h_0,h_0])}.
\end{cases}
\end{equation}
By $bc-ad>0$, it follows from  \cite[Theorem 7.1]{MM1996} that $\lim\limits_{t\to+\infty}(U(t),V(t))=(u^*,v^*)$. Moreover, the standard comparison principle yields that $(u(x,t),v(x,t))\preceq(U(t),V(t))$ for $x\in[g(t),h(t)]$ and $t>0$,  which leads to
\begin{equation}\label{2.9}
\limsup_{t\to+\infty}\left[\max\limits_{x\in[g(t),h(t)]}(u(x,t),v(x,t))\right]\preceq(u^*,v^*).
\end{equation}

We next analyze $\liminf\limits_{t\to+\infty}(u(x,t),v(x,t))$. Similar to the proof of  \cite[Theorem 3.6]{WD2021}, one can find a positive constant $L^{\ast}$ such that  spreading occurs when $h_{0}\geq L^{\ast}$. For any $l>L^{\ast}$, since $h_{\infty}=-g_{\infty}=+\infty$ and $g'(t)<0<h'(t)$ for $t>0$, it follows that there exists a constant $t_{l}>0$ such that
\begin{equation}\label{*1}
[-l,l]\subset[g(t),h(t)]\text{ for }t\ge t_l.
\end{equation}

Let $\Lambda_{0}$ be  the principal eigenvalue of the following cooperative eigenvalue problem:
\begin{equation*}
\begin{cases}
-d_{1}\omega_{xx}+a\omega-b\nu=\Lambda\omega, & -l<x<l,\\
-d_{2}\nu_{xx}-c\omega+d\nu=\Lambda\nu, & -l<x<l,\\
\omega(x)=\nu(x)=0, &x=\pm l.
\end{cases}
\end{equation*}
According to \cite[Lemma 3.3]{WD2021},  $\Lambda_{0}$ is negative and corresponds to a  positive eigenfunction $(\omega(x),\nu(x))$ given by $(\delta_1\cos(\pi x/(2l)),\cos(\pi x/(2l)))$ for some $\delta_1>0$.
Define
$$\underline{u}_{0}(x)=\tilde{\delta}\omega(x)~\textrm{and}~\underline{v}_{0}(x)=\tilde{\delta}\nu(x)\text{ for }x\in[-l,l],$$
where $\tilde{\delta}>0$ is a small constant to be determined. In view of $\mathrm{(A1)}$, we have $F(\tilde{\delta}\omega)=o_{\tilde{\delta}}(1)\tilde{\delta}\omega$ with $o_{\tilde{\delta}}(1)\to0^+$ as $\tilde{\delta}\to0^+$. Hence for $\tilde{\delta}>0$ sufficiently small, we have $o_{\tilde{\delta}}(1)+\Lambda_0\le0$, which leads to
\begin{equation}\label{Q0}
\begin{split}
&-d_{1}(\underline{u}_{0})_{xx}+a\underline{u}_{0}-b\underline{v}_{0}+F(\underline{u}_{0})
=\Lambda_{0}\tilde{\delta}\omega+o_{\tilde{\delta}}(1)\tilde{\delta}\omega\leq0\text{ for }x\in(-l,l),\\
&-d_{2}(\underline{v}_{0})_{xx}-c\underline{u}_{0}+d\underline{v}_{0}+G(\underline{v}_{0})
=\Lambda_{0}\tilde{\delta}\nu+o_{\tilde{\delta}}(1)\tilde{\delta}\nu\leq0\text{ for }x\in(-l,l).
\end{split}
\end{equation}
For any $l_{0}>l$,  \eqref{*1} yields that $[-l_0,l_0]\subset[g(t_{l_{0}}),h(t_{l_{0}})]$. Hence $(u(x,t_{l_{0}}),v(x,t_{l_{0}}))\succ(0,0)$ for $x\in[-l,l]$. Then we can choose $\tilde{\delta}>0$ smaller if necessary to obtain
\begin{equation}\label{Q00}
(\underline{u}_{0}(x),\underline{v}_{0}(x))=(\tilde{\delta}\omega(x),\tilde{\delta}\nu(x))\preceq(u(x,t_{l_{0}}),v(x,t_{l_{0}}))\text{ for }x\in[-l,l].
\end{equation}

We now consider the following auxiliary problem:
\begin{equation}\label{Q22}
\begin{cases}
(\underline{u}_{l})_{t}=d_{1}(\underline{u}_{l})_{xx}-a\underline{u}_{l}+b\underline{v}_{l}
-F(\underline{u}_{l}), & -l<x<l,~t>0,\\
(\underline{v}_{l})_{t}=d_{2}(\underline{v}_{l})_{xx}+c\underline{u}_{l}-d\underline{v}_{l}
-G(\underline{v}_{l}), & -l<x<l,~t>0,\\
\underline{u}_{l}(x,t)=\underline{v}_{l}(x,t)=0, & x=\pm l,~t>0,\\
\underline{u}_{l}(x,0)=\underline{u}_{0}(x),~\underline{v}_{l}(x,0)=\underline{v}_{0}(x), & -l\leq x\leq l.
\end{cases}
\end{equation}
By  \eqref{Q0}, we know that $(\underline{u}_0(x),\underline{v}_0(x))$ is a lower solution of system \eqref{Q22}. Moreover,  \eqref{*1}  and \eqref{Q00} imply that  $(u(x,t+t_{l_{0}}),v(x,t+t_{l_{0}}))$ is an upper solution  of system \eqref{Q22}. Thus,  one can derive from the standard parabolic comparison principle of cooperative systems that
\begin{equation}\label{2.11}
(u(x,t+t_{l_{0}}),v(x,t+t_{l_{0}}))\succeq(\underline{u}_{l}(x,t),\underline{v}_{l}(x,t))\text{ for }x\in[-l,l],~t>0.
\end{equation}

Furthermore, it follows from the theory of monotone dynamical systems (see, e.g., \cite{Smith1995}) that
\begin{equation}\label{Q3}
\lim_{t\to+\infty}(\underline{u}_{l}(x,t),\underline{v}_{l}(x,t))=
(\underline{u}_{l}^{\ast}(x),\underline{v}_{l}^{\ast}(x))\text{ uniformly in }[-l,l],
\end{equation}
where $(\underline{u}_{l}^{\ast}(x),\underline{v}_{l}^{\ast}(x))$ is the minimal solution of the following problem satisfying $(\underline{u}_{l}^{\ast}(x),\underline{v}_{l}^{\ast}(x))\succeq(\underline{u}_0(x),\underline{v}_0(x))$ for $x\in[-l,l]$:
\begin{equation*}\label{2.12}
\begin{cases}
d_{1}(\underline{u}_{l}^{\ast})_{xx}-a\underline{u}_{l}^{\ast}+b\underline{v}_{l}^{\ast}
-F(\underline{u}_{l}^{\ast})=0, & -l<x<l,\\
d_{2}(\underline{v}_{l}^{\ast})_{xx}+c\underline{u}_{l}^{\ast}-d\underline{v}_{l}^{\ast}
-G(\underline{v}_{l}^{\ast})=0, & -l<x<l,\\
\underline{u}_{l}^{\ast}(x)=\underline{u}_{l}^{\ast}(x)=0, & x=\pm l.
\end{cases}
\end{equation*}

For any $l_{1}\geq l_{2}>L^{\ast}$, by comparing the boundary conditions in system  \eqref{Q22}  with $l=l_1$ and $l=l_2$, we derive from the standard parabolic comparison principle of cooperative systems that
$$(\underline{u}_{l_1}(x,t),\underline{v}_{l_1}(x,t))\succeq(\underline{u}_{l_2}(x,t),\underline{v}_{l_2}(x,t))\text{ for }x\in[-l_2,l_2],~t>0.$$
Thus, letting $t\to+\infty$ and using \eqref{Q3}, we obtain  $(\underline{u}_{l_{1}}^*(x),\underline{v}_{l_{1}}^*(x))\succeq(\underline{u}_{l_{2}}^*(x),\underline{v}_{l_{2}}^*(x))$ for $x\in[-l_{2},l_{2}]$. This means that the pair $(\underline{u}_l^*(x),\underline{v}_l^*(x))$  is non-decreasing with respect to $l$.

By the classical elliptic regularity theory and  diagonal procedure, we know that
 \begin{equation}\label{eq02}
\lim\limits_{l\to+\infty}(\underline{u}_{l}^{\ast}(x),\underline{v}_{l}^{\ast}(x))=(\underline{u}_{\infty}^{\ast}(x),\underline{v}_{\infty}^{\ast}(x))\text{ in }\big[C^{2}_{\mathrm{loc}}(\mathbb{R})\big]^{2},
 \end{equation}
where $(\underline{u}_{\infty}^{\ast},\underline{v}_{\infty}^{\ast})$ satisfies
\begin{equation}\label{2.13}
\begin{cases}
d_{1}(\underline{u}_{\infty}^{\ast})_{xx}-a\underline{u}_{\infty}^{\ast}+b\underline{v}_{\infty}^{\ast}
-F(\underline{u}_{\infty}^{\ast})=0, & x\in\mathbb{R},\\
d_{2}(\underline{v}_{\infty}^{\ast})_{xx}+c\underline{u}_{\infty}^{\ast}-d\underline{v}_{\infty}^{\ast}
-G(\underline{v}_{\infty}^{\ast})=0, & x\in\mathbb{R}.
\end{cases}
\end{equation}

We claim that
 \begin{equation}\label{eq03}
 (\underline{u}_{\infty}^{\ast}(x),\underline{v}_{\infty}^{\ast}(x))\equiv(u^*,v^*)\text{ for }x\in\mathbb{R}.
 \end{equation}
Indeed, for any $\hat{x}\in(-\frac{l}{2},\frac{l}{2})$, denote $l_{1}:=|\hat{x}|$. Then
$[-l+2l_1,l-2l_1]\subset[-l+l_1-\hat{x},l-l_1-\hat{x}]\subset[-l,l]$.
Hence it follows from the monotonicity of $(\underline{u}_l^*(x),\underline{v}_l^*(x))$ in $l$ that
$$(\underline{u}_{l-2l_{1}}^{\ast}(x),\underline{v}_{l-2l_{1}}^{\ast}(x))\preceq
(\underline{u}_{l-l_1}^{\ast}(x+\hat{x}),\underline{v}_{l-l_1}^{\ast}(x+\hat{x}))\preceq
(\underline{u}_{l}^{\ast}(x),\underline{v}_{l}^{\ast}(x)).$$
Letting $l\to+\infty$, we obtain
$$(\underline{u}_{\infty}^{\ast}(x),\underline{v}_{\infty}^{\ast}(x))\preceq
(\underline{u}_{\infty}^{\ast}(x+\hat{x}),\underline{v}_{\infty}^{\ast}(x+\hat{x}))\preceq
(\underline{u}_{\infty}^{\ast}(x),\underline{v}_{\infty}^{\ast}(x))\text{ for }x\in\mathbb{R},$$
which leads to $(\underline{u}_{\infty}^{\ast}(\hat{x}),\underline{v}_{\infty}^{\ast}(\hat{x}))=
(\underline{u}_{\infty}^{\ast}(0),\underline{v}_{\infty}^{\ast}(0))$.
The arbitrariness of $\hat{x}$ yields
$$(\underline{u}_{\infty}^{\ast}(x),\underline{v}_{\infty}^{\ast}(x))\equiv
(\underline{u}_{\infty}^{\ast}(0),\underline{v}_{\infty}^{\ast}(0))\text{ for }x\in\mathbb{R}.$$
By $bc-ad>0$,  $(u^*,v^*)$ is the unique positive equilibrium solution of \eqref{2.13}. Thus, \eqref{eq03} holds.

From \eqref{eq02} and \eqref{eq03}, for any $\epsilon>0$ and $L_*>0$, there exists a constant $l>\max\{L_*,L^{\ast}\}$ such that
\begin{equation*}\label{2.14}
(\underline{u}_{l}^{\ast}(x),\underline{v}_{l}^{\ast}(x))\succeq(u^*-\epsilon,u^*-\epsilon)\text{ for }x\in[-L_*,L_*].
\end{equation*}
Combining  with \eqref{2.11} and \eqref{Q3}, we obtain
$$\liminf_{t\to+\infty}(u(x,t),v(x,t))\succeq(u^*-\epsilon,u^*-\epsilon)
~\textrm{uniformly~in}~[-L_*,L_*].$$
By the arbitrariness of $\epsilon$ and $L_*$, we have
$$\liminf_{t\to+\infty}(u(x,t),v(x,t))\succeq(u^*,v^*)
~\textrm{uniformly~in~any~compact~subset~of}~\mathbb{R}.$$
This, together with \eqref{2.9}, yields the desired result. The proof is complete.
\end{proof}
\bigskip
\textbf{Proof of Theorem \ref{theo1.1}.}
It follows from Lemmas  \ref{lem2.5}  and  \ref{lem2.6}  that Theorem \ref{theo1.1} holds.
 $\hfill\square$

\section{Asymptotic spreading speed}\label{section 4}
\setcounter{equation}{0}

This section, which examines the dynamics of the free boundaries $h(t)$ and $g(t)$ in the spreading scenario, is divided into two subsections. In Subsection  \ref{subsection 4.1}, we complete the proofs of Proposition  \ref{theo1.2} and Theorem \ref{theo1.3}. Subsection \ref{subsection 4.2} is dedicated to showing the boundedness of $h(t)-s_{\mu,\rho}t$ and $g(t)+s_{\mu,\rho}t$  for $t\ge0$, which plays  a crucial role in the proof of   Theorem \ref{theo1.4}.

\subsection{Proofs of  Proposition  \ref{theo1.2} and Theorem \ref{theo1.3}}\label{subsection 4.1}
We first prove  Proposition  \ref{theo1.2}. For any $s\ge0$, define
$$p_{1,s}(\lambda):=d_{1}\lambda^{2}-s\lambda-a,\quad p_{2,s}(\lambda):=d_{2}\lambda^{2}-s\lambda-d$$
and
$$P_{s}(\lambda):=p_{1,s}(\lambda)p_{2,s}(\lambda)-bc=(d_{1}\lambda^{2}-s\lambda-a)(d_{2}\lambda^{2}-s\lambda-d)-bc.$$
Clearly,    the polynomial $p_{i,s}(\lambda)$ has two roots  $\lambda_i^-$ and $\lambda_i^+$ satisfying
$$\lambda_{i}^-<0<\lambda_i^+\text{ and }p_{i,s}(\lambda)<0\text{ in }(\lambda_{i}^{-},\lambda_{i}^{+})\text{ for }i=1,2.$$
By the definition  of $P_s(\lambda)$, we obtain
\begin{equation}\label{add equation 02}
P_s(0)=ad-bc<0,~P_s(\pm\infty)=+\infty,~P_s(\lambda_i^{\pm})=-bc<0 \text{ for }i=1,2.
\end{equation}
A direct calculation yields that $\lim\limits_{s\to+\infty}P_s(1/\sqrt{s})=+\infty$.  Hence for all large $s$, there holds
\begin{equation*}
P_s(1/\sqrt{s})>0\text{ and }0<1/\sqrt{s}<\lambda_i^+\text{ for }i=1,2.
\end{equation*}
Combining   with \eqref{add equation 02}, we conclude that $P_s(\lambda)$ has four real roots for all large $s$. Moreover, it is easy to verify that $P_0(\lambda)$ has a pair of conjugate pure imaginary roots. Therefore,
$$s^{\ast}:=\inf\left\{s_{0}>0:\text{ All roots  of }P_s(\lambda)\text{ are real for }s\geq s_{0}\right\}$$
is well-defined, and $s^{\ast}>0$.

\begin{lem}\label{lemmap1}
 Assume that  $\mathrm{(A1)}$-$\mathrm{(A3)}$ hold and $bc-ad>0$. Then for every $s\in[0, s^{\ast})$,
system \eqref{1.5} has a unique strictly increasing solution $(\phi_s,\psi_s)\in C^2([0,+\infty))\times C^2([0,+\infty))$, and system  \eqref{1.5}  has no such solution when $s\geq s^{\ast}$. Moreover, $\lim\limits_{s\to(s^*)^-}(\phi_s,\psi_s)=(0,0)$ in $C^2_{\mathrm{loc}}([0,+\infty))\times C^2_{\mathrm{loc}}([0,+\infty))$.
\end{lem}
\begin{proof}
The proof follows a similar approach to that of \cite[Theorem 3.2]{WND2019} and is  therefore omitted here.
\end{proof}
\begin{lem}\label{lemmap2}
Assume that  $\mathrm{(A1)}$-$\mathrm{(A3)}$ hold and $bc-ad>0$. Then for any $\mu>0$ and $\rho>0$, there exists a unique $s_{\mu,\rho}\in(0,s^{\ast})$ such that \eqref{s mu} holds.
\end{lem}
\begin{proof}
Let $(\phi_{s},\psi_{s})$ be the strictly increasing solution of problem \eqref{1.5} with $s\in[0,s^{\ast})$. We  claim that
\begin{equation}\label{2.22}
\phi_{s_{1}}(\xi)>\phi_{s_{2}}(\xi),~\psi_{s_{1}}(\xi)>\psi_{s_{2}}(\xi)\text{ in }(0,+\infty)\text{ when }0\le s_{1}<s_{2}<s^{*}.
\end{equation}

Indeed, in light of $0\leq s_{1}<s_{2}$ and $\phi'_{s_{2}}>0, \psi'_{s_{2}}>0$ in $\mathbb{R}^{+}$, we have
\begin{equation*}
\begin{cases}
d_1 \phi''_{s_{2}}-s_{1}\phi_{s_{2}}'-a\phi_{s_{2}}+b\psi_{s_{2}}-F(\phi_{s_{2}})=(s_2-s_1)\phi'_{s_2}>0, & 0<\xi<+\infty,\\
d_2 \psi''_{s_{2}}-s_{1}\psi_{s_{2}}'+c\phi_{s_{2}}-d\psi_{s_{2}}-G(\psi_{s_{2}})=(s_2-s_1)\psi'_{s_2}>0, & 0<\xi<+\infty
\end{cases}
\end{equation*}
with $(\phi_{s_{i}}(0),\psi_{s_{i}}(0))=(0,0)$ and $(\phi_{s_{i}}(+\infty),\psi_{s_{i}}(+\infty))=(u^*,v^*)$ for $i=1,2$. Therefore,
it follows from the uniqueness of the monotone solution  that
 $$(\phi_{s_{1}}(\xi),\psi_{s_{1}}(\xi))\succeq(\phi_{s_{2}}(\xi),\psi_{s_{2}}(\xi))\text{ for }0<\xi<+\infty.$$
Furthermore, by the strong maximum principle, the strict inequality \eqref{2.22} holds.

For any given $\mu>0$ and $\rho>0$, we define
$$f(s):=\mu[\psi_{s}'(0)+\rho\phi_{s}'(0)]-s\text{ for }s\in[0,s^*).$$
From \eqref{2.22}, both $\phi_s'(0)$ and $\psi_s'(0)$  are  non-increasing with respect to $s\in[0,s^*)$. Thus,
$f(s)$ is continuous and strictly decreasing with respect to $s\in[0,s^{\ast})$.
It follows from Lemma \ref{lemmap1} that $\phi_0'(0)>0$, $\psi_0'(0)>0$ and $\phi_{s^*}'(0)=\psi_{s^*}'(0)=0$, which implies that $f(0)=\mu[\psi_{0}'(0)+\rho\phi_{0}'(0)]>0$ and $\lim\limits_{s\to ({s^{\ast}})^{-}}f(s)=-s^{\ast}<0$.
Therefore, for any $\mu>0$ and $\rho>0$,  there exists a unique $s_{\mu,\rho}\in(0,s^{\ast})$ such that $f(s_{\mu,\rho})=0$, namely, \eqref{s mu} is valid. This completes the proof.
\end{proof}

\textbf{Proof of  Proposition  \ref{theo1.2}.}
 It follows from  Lemmas \ref{lemmap1} and \ref{lemmap2} that  Proposition  \ref{theo1.2} holds.
 $\hfill\square$

 \bigskip

Next, we derive an exponential estimate for the semi-wave solution of  system \eqref{1.5}.
\begin{lem}\label{lem2.4}
 Assume that  $\mathrm{(A1)}$-$\mathrm{(A3)}$ hold and $bc-ad>0$. Let $(\phi(\xi),\psi(\xi))$ be a monotone solution of system \eqref{1.5}. Then there exist  constants $\hat{\mu}_{1}, p, q>0$ such that as $\xi \to+\infty$,
\begin{equation*}
\begin{cases}
(\phi(\xi),\psi(\xi))=(u^*,v^*)-e^{-\hat{\mu}_{1}\xi}(p+o(1),q+o(1)),\\
(\phi'(\xi),\psi'(\xi))=(O(e^{-\hat{\mu}_{1}\xi}),O(e^{-\hat{\mu}_{1}\xi})).
\end{cases}\label{2.2}
\end{equation*}
\end{lem}
\begin{proof}
A simple calculation shows that $(u^\ast, 0, v^\ast, 0)$ is a saddle point of the ODE system satisfied by $(\phi,\phi',\psi,\psi')$. Then  the standard theory of stable and unstable manifolds implies that both $u^\ast-\phi$ and $v^\ast-\psi$ converge to zero exponentially
as $\xi\rightarrow+\infty$. For the remainder of the proof, one can refer to \cite[Lemma 2.3]{DWZ2017} for details.
\end{proof}
\bigskip\textbf{Proof of Theorem \ref{theo1.3}.}
Inspired by \cite[Theorem 3.15]{WND2019}, we now prove that
$\lim\limits_{t \to +\infty}h(t)/t=s_{\mu,\rho}$. The same argument applied to $(u(-x,t),v(-x,t),-h(t),-g(t))$  yields the other limit.
For clarity,  the proof is divided into two steps.

\textbf{Step 1.} We first prove that
\begin{equation}\label{limsup of h(t)/t}
\limsup_{t \to +\infty}\frac{h(t)}{t}\leq s_{\mu,\rho}.
\end{equation}

Indeed, for any $\tau\in(0,\min\{a,d\})$, let $a_{\tau}=a-\tau$ and $d_{\tau}=d-\tau$. Then  $bc-a_{\tau}d_{\tau}>bc-ad>0$. Hence system \eqref{1.11} with $(a,d)=(a_{\tau},d_{\tau})$ has a unique positive equilibrium, denoted by $(u^*_{\tau},v^*_{\tau})$. Similar to the proof of Lemma \ref{lem2.6}, one can show that \eqref{2.9} holds with $(u^*,v^*)=(u^*_{\tau},v^*_{\tau})$, that is,
\begin{equation*}
\limsup_{t\to+\infty}\left[\max\limits_{x\in[g(t),h(t)]}(u(x,t),v(x,t))\right]\preceq(u^*_{\tau},v^*_{\tau}).
\end{equation*}
Thus, for given  $\epsilon_1>0$ small, there exists $T_{1}>0$ large enough such  that
\begin{equation}\label{1.13}
(u(x,t),v(x,t))\preceq(u^*_{\tau}+\epsilon_1,v^*_{\tau}+\epsilon_1)\text{ for }x\in[g(t),h(t)],~t\geq T_{1}.
\end{equation}
Similarly, system \eqref{1.11} with $(a,d)=(a_{2\tau},d_{2\tau})$ has a unique positive equilibrium, denoted by  $(u^*_{2\tau},v^*_{2\tau})$. Clearly,  there holds $(u^*_{\tau},v^*_{\tau})\prec(u^*_{2\tau},v^*_{2\tau})$.

Consider the following auxiliary problem:
\begin{equation}
\begin{cases}
d_1 \phi_{\tau}''-s_{\mu,\rho}^{\tau}\phi_{\tau}'-(a-2\tau)\phi_{\tau}+b\psi_{\tau}-F(\phi_{\tau})=0,~ \phi_{\tau}'>0, &0<\xi<+\infty,\\
d_2 \psi_{\tau}''-s_{\mu,\rho}^{\tau}\psi_{\tau}'+c\phi_{\tau}-(d-2\tau)\psi_{\tau}-G(\psi_{\tau})=0,~ \psi_{\tau}'>0, &0<\xi<+\infty,\\
(\phi_{\tau}(0),\psi_{\tau}(0))=(0,0),~ (\phi_{\tau}(+\infty),\psi_{\tau}(+\infty))=(u^*_{2\tau},v^*_{2\tau}).
\end{cases}
\label{1.12}
\end{equation}
By Proposition  \ref{theo1.2}, there exists a unique $s_{\mu,\rho}^{\tau}>0$ such that problem \eqref{1.12} has a unique  strictly increasing solution $(\phi_{\tau},\psi_{\tau})$ satisfying
\begin{equation}\label{equation 2}
\mu[\psi_{\tau}'(0)+\rho \phi_{\tau}'(0)]=s_{\mu,\rho}^{\tau}\text{ and }\lim_{\tau\to0^{+}}s_{\mu,\rho}^{\tau}=s_{\mu,\rho}.
\end{equation}
In view of $(\phi_{\tau}(+\infty),\psi_{\tau}(+\infty))=(u^*_{2\tau},v^*_{2\tau})\succ(u^*_{\tau},v^*_{\tau})$, we can choose $\xi_{0}>0$ large enough and $\epsilon_1>0$ smaller if necessary to obtain that
\begin{equation*}
(\phi_{\tau}(\xi_{0}),\psi_{\tau}(\xi_{0}))\succ(u^*_{\tau}+\epsilon_1,v^*_{\tau}+\epsilon_1).
\end{equation*}
Combining with \eqref{1.13}, we have
\begin{equation}\label{1.14}
(\phi_{\tau}(\xi_{0}),\psi_{\tau}(\xi_{0}))\succ(u(x,t),v(x,t))\text{ for }x\in[g(t),h(t)],~t\geq T_{1}.
\end{equation}

Define
\begin{equation*}
\begin{split}
&\tilde{h}(t)=s_{\mu,\rho}^{\tau}(t-T_{1})+\xi_{0}+h(T_{1})~\text{ for }~t\geq T_{1},\\
&\tilde{u}(x,t)=\phi_{\tau}(\tilde{h}(t)-x)~\text{ for }~0\leq x\leq \tilde{h}(t),~t\geq T_{1},\\
&\tilde{v}(x,t)=\psi_{\tau}(\tilde{h}(t)-x)~\text{ for }~0\leq x\leq \tilde{h}(t),~t\geq T_{1}.
\end{split}
\end{equation*}
Then $\tilde{h}(T_1)>h(T_{1})$ and  $\tilde{u}(\tilde{h}(t),t)=\tilde{v}(\tilde{h}(t),t)=0$ for $t\ge T_1$. Moreover, \eqref{equation 2} yields
$$\tilde{h}'(t)=s_{\mu,\rho}^{\tau}=\mu [\psi_{\tau}'(0)+\rho \phi_{\tau}'(0)]=-\mu[\tilde{v}_{x}(\tilde{h}(t),t)+\rho\tilde{u}_{x}(\tilde{h}(t),t)]
\text{ for }t\geq T_{1}.$$
In view of \eqref{1.14}, it follows from the strict monotonicity of $\phi_{\tau}$ and $\psi_{\tau}$ that
\begin{equation*}
\begin{split}
&\tilde{u}(x,T_{1})=\phi_{\tau}(\xi_{0}+h(T_{1})-x)\geq \phi_{\tau}(\xi_{0})>u(x,T_{1})\text{ for }x\in[0,h(T_{1})],\\
&\tilde{v}(x,T_{1})=\psi_{\tau}(\xi_{0}+h(T_{1})-x)\geq \psi_{\tau}(\xi_{0})>v(x,T_{1})\text{ for }x\in[0,h(T_{1})],\\
&\tilde{u}(0,t)=\phi_{\tau}(\tilde{h}(t))\geq \phi_{\tau}(\xi_{0})>u(0,t)\text{ for }t\ge T_1,\\
 &\tilde{v}(0,t)=\psi_{\tau}(\tilde{h}(t))\geq \psi_{\tau}(\xi_{0})>v(0,t)\text{ for }t\ge T_1.
\end{split}
\end{equation*}
Finally, some direct calculations show that
for  $x\in(0,\tilde{h}(t))$ and $t\ge T_1$,
\begin{equation*}
\begin{split}
&\tilde{u}_{t}-d_{1}\tilde{u}_{xx}=s_{\mu,\rho}^{\tau}\phi_{\tau}'-d_{1}\phi_{\tau}''
=-(a-2\tau)\phi_{\tau}+b\psi_{\tau}-F(\phi_{\tau})
\geq -a\tilde{u}+b\tilde{v}-F(\tilde{u}),\\
&\tilde{v}_{t}-d_{2}\tilde{v}_{xx}=s_{\mu,\rho}^{\tau}\psi_{\tau}'-d_{2}\psi_{\tau}''
=c\phi_{\tau}-(d-2\tau)\psi_{\tau}-G(\psi_{\tau})
\geq  c\tilde{u}-d\tilde{v}-G(\tilde{v}).
\end{split}
\end{equation*}
 Therefore, it follows from Lemma \ref{lem2.1} that
$$h(t)\leq\tilde{h}(t)\text{ for }t\ge T_1,$$
 which leads to  $\limsup\limits_{t\to+\infty}h(t)/t\leq s_{\mu,\rho}^{\tau}$.
 Letting  $\tau\to0^{+}$ and using \eqref{equation 2}, we obtain \eqref{limsup of h(t)/t}.

\textbf{Step 2.} We now prove that
\begin{equation}\label{liminf of h(t)/t}
\liminf_{t \to +\infty}\frac{h(t)}{t}\geq s_{\mu,\rho}.
\end{equation}

Indeed, for any $\epsilon>0$, let $a_{\epsilon}=a+\epsilon$ and $d_{\epsilon}=d+\epsilon$. In light of $bc-ad>0$, we have   $bc-a_{\epsilon}d_{\epsilon}>0$ for small $\epsilon$. Hence system \eqref{1.11} with $(a,d)=(a_{\epsilon},d_{\epsilon})$ has a unique positive equilibrium, denoted by $(u_{\epsilon}^*,v_{\epsilon}^*)$. Clearly, we have $(u^*,v^*)\succ(u_{\epsilon}^*,v_{\epsilon}^*)$. Since spreading occurs, it follows that $\lim\limits_{t\to+\infty}h(t)=-\lim\limits_{t\to+\infty}g(t)=+\infty$ and $\lim\limits_{t\to+\infty}(u(x,t),v(x,t))=(u^*,v^*)$ locally uniformly in $\mathbb{R}$. Thus, there exist constants $T_2>0$  and $K\in(0,h(T_2))$ such that
\begin{equation}\label{equation 333}
h(t)>K\text{ and }(u(x,t),v(x,t))\succ(u_{\epsilon}^*,v_{\epsilon}^*)\text{ for }x\in[0,K],~t\ge T_2.
\end{equation}

Consider the following auxiliary problem:
\begin{equation}
\begin{cases}
d_1 \phi_{\epsilon}''-s_{\mu,\rho}^{\epsilon}\phi_{\epsilon}'-(a+\epsilon)\phi_{\epsilon}+b\psi_{\epsilon}-F(\phi_{\epsilon})=0,~ \phi_{\epsilon}'>0, &0<\xi<+\infty,\\
d_2 \psi_{\epsilon}''-s_{\mu,\rho}^{\epsilon}\psi_{\epsilon}'+c\phi_{\epsilon}-(d+\epsilon)\psi_{\epsilon}-G(\psi_{\epsilon})=0,~ \psi_{\epsilon}'>0, &0<\xi<+\infty,\\
(\phi_{\epsilon}(0),\psi_{\epsilon}(0))=(0,0),~ (\phi_{\epsilon}(+\infty),\psi_{\epsilon}(+\infty))=(u_{\epsilon}^*,v_{\epsilon}^*).
\end{cases}
\label{1.15}
\end{equation}
By Proposition $\ref{theo1.2}$, there exists a unique $s_{\mu,\rho}^{\epsilon}>0$ such that problem \eqref{1.15} has a  strictly increasing solution $(\phi_{\epsilon},\psi_{\epsilon})$ satisfying
\begin{equation}\label{equation 3}
\mu[\psi_{\epsilon}'(0)+\rho \phi_{\epsilon}'(0)]=s_{\mu,\rho}^{\epsilon}\text{ and }\lim_{\epsilon\to0^{+}}s_{\mu,\rho}^{\epsilon}=s_{\mu,\rho}.
\end{equation}

Define
\begin{equation*}
\begin{split}
&\check{h}(t)=s_{\mu,\rho}^{\epsilon}(t-T_{2})+K~\text{ for }~t\geq T_{2},\\
&\check{u}(x,t)=\phi_{\epsilon}(\check{h}(t)-x)~\text{ for }~0\leq x\leq \check{h}(t),~t\geq T_{2},\\
&\check{v}(x,t)=\psi_{\epsilon}(\check{h}(t)-x)~\text{ for }~0\leq x\leq \check{h}(t),~t\geq T_{2}.
\end{split}
\end{equation*}
Then $\check{h}(T_{2})=K<h(T_{2})$ and $\check{u}(\check{h}(t),t)=\check{v}(\check{h}(t),t)=0$ for $t\ge T_2$. Due to \eqref{equation 3}, we get
$$\check{h}'(t)=s_{\mu,\rho}^{\epsilon}=\mu[\psi_{\epsilon}'(0)+\rho \phi_{\epsilon}'(0)]=-\mu[\check{v}_{x}(\check{h}(t),t)+\rho \check{u}_{x}(\check{h}(t),t)]\text{ for }t\geq T_{2}.$$
A direct calculation yields that
for $x\in(0,\check{h}(t))$ and $t\ge T_2$,
\begin{equation*}
\begin{split}
&\check{u}_{t}-d_{1}\check{u}_{xx}
= s_{\mu,\rho}^{\epsilon}\phi_{\epsilon}'-d_{1}\phi_{\epsilon}''
=-(a+\epsilon)\phi_{\epsilon}+b\psi_{\epsilon}-F(\phi_{\epsilon})
\leq -a\check{u}+b\check{v}-F(\check{u}),\\
&\check{v}_{t}-d_{2}\check{v}_{xx}
= s_{\mu,\rho}^{\epsilon}\psi_{\epsilon}'-d_{2}\psi_{\epsilon}''
=c\phi_{\epsilon}-(d+\epsilon)\psi_{\epsilon}-G(\psi_{\epsilon})
\leq c\check{u}-d\check{v}-G(\check{v}).
\end{split}
\end{equation*}
Moreover, it follows from \eqref{equation 333}  that
\begin{equation*}
\begin{split}
&\check{u}(x,T_{2})=\phi_{\epsilon}(h(T_{2})-x)< u_{\epsilon}^*<u(x,T_{2})\text{ for }x\in[0,K]=[0,\check{h}(T_2)],\\
&\check{v}(x,T_{2})=\psi_{\epsilon}(h(T_{2})-x)< v_{\epsilon}^*<v(x,T_{2})\text{ for }x\in[0,K]=[0,\check{h}(T_2)],\\
&\check{u}(0,t)=\phi_{\epsilon}(\check{h}(t))< u_{\epsilon}^*< u(0,t)\text{ for }t\ge T_2,\\
&\check{v}(0,t)=\psi_{\epsilon}(\check{h}(t))< v_{\epsilon}^*<v(0,t)\text{ for }t\ge T_2.
\end{split}
\end{equation*}
Hence by virtue of Remark \ref{remark2.2}(i), we obtain
$$h(t)\ge\check{h}(t)\text{ for }t\geq T_{2},$$
which implies that $\liminf\limits_{t\to+\infty}h(t)/t\geq s_{\mu,\rho}^{\epsilon}$. Letting $\epsilon\to 0^+$ and using \eqref{equation 3}, we get \eqref{liminf of h(t)/t}.

Therefore, it follows from \eqref{limsup of h(t)/t} and \eqref{liminf of h(t)/t} that $\lim\limits_{t \to +\infty}h(t)/t=s_{\mu,\rho}$. The proof is complete.
 $\hfill\square$

\subsection{Bounds for $g(t)+s_{\mu,\rho}t$ and $h(t)-s_{\mu,\rho}t$} \label{subsection 4.2}
This subsection is devoted to showing that both $g(t)+s_{\mu,\rho}t$ and $h(t)-s_{\mu,\rho}t$ are bounded functions for $t\ge0$. More precisely, our goal is to achieve  the following result.
\begin{prop}\label{prop3.1}
Assume that $\mathrm{(H)}$ holds and let $(u,v,g,h)$ be the unique solution of system \eqref{1.2}. If spreading occurs, then there exists a  constant $C>0$ such that
$$|{g(t)+s_{\mu,\rho}t}|\le C\text{ and }|{h(t)-s_{\mu,\rho}t}|\leq C\text{ for all }t\ge 0.$$
\end{prop}

Motivated by the work in \cite{WND}, we prove this result by constructing appropriate upper and lower solutions. Firstly, we construct an upper solution $(\overline{u},\overline{v},\overline{g},\overline{h})$ of system \eqref{1.2} as follows:
\begin{equation}
\begin{split}
&\overline{h}(t)=s_{\mu,\rho}(t-T^{\ast})+\sigma(1-e^{-\delta(t-T^{\ast})})+h(T^{\ast})+X_0,\quad \overline{g}(t)=g(t),\\
&\overline{u}(x,t)=(1+K_{1}e^{-\delta(t-T^{\ast})})\phi_{s_{\mu,\rho}}(\overline{h}(t)-x),\\
&\overline{v}(x,t)=(1+K_{1}e^{-\delta(t-T^{\ast})})\psi_{s_{\mu,\rho}}(\overline{h}(t)-x),
\end{split}\label{3.1}
\end{equation}
where $(\phi_{s_{\mu,\rho}},\psi_{s_{\mu,\rho}})$ is  the strictly increasing solution of system \eqref{1.5} with $s=s_{\mu,\rho}$,   and $X_0, T^{\ast}, K_{1}, \delta, \sigma$ are positive constants to be determined.
\begin{lem}\label{lem3.2}
Assume that $\mathrm{(H)}$ holds and let $(u,v,g,h)$ be the unique solution of system \eqref{1.2}. If spreading occurs,  then  there  exist positive constants $X_0, T^{\ast},  K_{1}, \delta$, and $\sigma$  such that
$$h(t)\le\overline{h}(t),~(u(x,t),v(x,t))\preceq(\overline{u}(x,t),\overline{v}(x,t))\text{ for }x\in[g(t),h(t)],~t>T^*.$$
\end{lem}
\begin{proof}
We claim that $(\overline{u},\overline{v},\overline{g},\overline{h})$ is an upper solution of system \eqref{1.2} for $t>T^{\ast}$ by taking appropriate parameters  $X_0, T^{\ast}, K_{1}, \delta$ and $\sigma$, that is,
\begin{align}
&\overline{u}_{t}\geq d_1 \overline{u}_{xx}-a\overline{u}+b\overline{v}-F(\overline{u}),&&\overline{g}(t)< x<\overline{h}(t),~\ t>T^{\ast},\label{3.2}\\
&\overline{v}_{t}\geq d_2 \overline{v}_{xx}+c\overline{u}-d\overline{v}-G(\overline{v}),&&\overline{g}(t)< x<\overline{h}(t),~\ t>T^{\ast},\label{3.3}\\
&\overline{u}(x,t)\geq u(x,t),~\overline{v}(x,t)\geq v(x,t),&&x=\overline{g}(t),~t>T^{\ast},\label{3.4}\\
&\overline{u}(x,t)=\overline{v}(x,t)=0,&&x=\overline{h}(t),~t>T^{\ast},\label{3.5}\\
&\overline{h}(T^{\ast})\geq h(T^{\ast}),~\overline{h}'(t)\geq -\mu [\overline{v}_{x}(\overline{h}(t),t)+\rho\overline{u}_{x}(\overline{h}(t),t)],&&t>T^{\ast},\label{3.6}\\
&\overline{u}(x,T^{\ast})\geq {u}(x,T^{\ast}),~\overline{v}(x,T^{\ast})\geq {v}(x,T^{\ast}),&&\overline{g}(T^{\ast})\leq x\leq h(T^{\ast}).\label{3.7}
\end{align}

Firstly, the definitions  of $\overline{u}$, $\overline{v}$, and $\overline{g}$ imply that for $t>T^*$,
$$(\overline{u}(\overline{g}(t),t),\overline{v}(\overline{g}(t),t))\succeq(0,0)=(u(\overline{g}(t),t),v(\overline{g}(t),t))
\text{ and }\overline{u}(\overline{h}(t),t)=
\overline{v}(\overline{h}(t),t)=0,$$
 which yields \eqref{3.4} and \eqref{3.5}.

Clearly,  $\overline{h}(T^{\ast})\geq h(T^{\ast})$ and $\overline{h}'(t)=s_{\mu,\rho}+\sigma\delta e^{-\delta(t-T^{\ast})}$ for $t>T^*$. Moreover, \eqref{s mu} gives rise to
$$-\mu[\overline{v}_{x}(\overline{h}(t),t)+\rho\overline{u}_{x}(\overline{h}(t),t)]=
s_{\mu,\rho}(1+K_{1}e^{-\delta(t-T^{\ast})})\text{ for }t>T^*.$$
Hence \eqref{3.6} holds provided that
\begin{equation}\label{3.8}
\sigma\delta\ge s_{\mu,\rho}K_{1}.
\end{equation}

Theorem \ref{lem2.0} suggests that there exists  $C_0\ge\max\{u^*,v^*\}$ such that  $u(x,t),v(x,t)\le C_0$  for $x\in[g(t),h(t)]$ and $t>0$. For such $C_0$, one can  take $K_{1}$ large enough to obtain
\begin{equation}\label{3.9}
(1+K_{1})\phi_{s_{\mu,\rho}}(X_0)\geq{C_{0}}\text{ and }(1+K_{1})\psi_{s_{\mu,\rho}}(X_0)\geq{C_{0}}.
\end{equation}
Combining \eqref{3.9} with the strict monotonicity of  $\phi_{s_{\mu,\rho}}$, it follows that for $x \in [\overline{g}(T^{\ast}),h(T^{\ast})]$,
\begin{equation*}
\overline{u}(x,T^{\ast})=(1+K_{1})\phi_{s_{\mu,\rho}}(h(T^{\ast})+X_0-x)
\geq (1+K_{1})\phi_{s_{\mu,\rho}}(X_0)\ge C_0
\geq u(x,T^{\ast}).
\end{equation*}
A similar argument yields that $\overline{v}(x,T^{\ast})\geq v(x,T^{\ast})$ for $x \in [\overline{g}(T^{\ast}),h(T^{\ast})]$, which proves \eqref{3.7}.

It remains to prove \eqref{3.2} and \eqref{3.3}.
 For convenience, we set $z=\overline{h}(t)-x$ for $x\in(g(t),\overline{h}(t))$ and $t>T^*$.
  Then
  $$z\in(0,\overline{h}(t)-g(t))\text{ with }\lim\limits_{t\to+\infty}\big(\overline{h}(t)-g(t)\big)=+\infty.$$
A direct calculation yields that
\begin{align*}
&\overline{u}_{t}=-\delta K_{1}e^{-\delta(t-T^{\ast})}\phi_{s_{\mu,\rho}}(z)+(1+K_{1}e^{-\delta(t-T^{\ast})})(s_{\mu,\rho}+\sigma\delta
e^{-\delta(t-T^{\ast})})\phi_{s_{\mu,\rho}}'(z),\\
&\overline{v}_{t}=-\delta K_{1}e^{-\delta(t-T^{\ast})}\psi_{s_{\mu,\rho}}(z)+(1+K_{1}e^{-\delta(t-T^{\ast})})(s_{\mu,\rho}+\sigma\delta
e^{-\delta(t-T^{\ast})})\psi_{s_{\mu,\rho}}'(z),\\
&\overline{u}_{xx}=(1+K_{1}e^{-\delta(t-T^{\ast})})\phi_{s_{\mu,\rho}}''(z),\quad
\overline{v}_{xx}=(1+K_{1}e^{-\delta(t-T^{\ast})})\psi_{s_{\mu,\rho}}''(z).
\end{align*}
By the equation of $\phi_{s_{\mu,\rho}}$ (see \eqref{1.5}), we have
\begin{equation}\label{ad ee}
\begin{split}
&\ \overline{u}_{t}-d_1 \overline{u}_{xx}+a\overline{u}-b\overline{v}+F(\overline{u})\\
=&-\delta K_{1}e^{-\delta(t-T^{\ast})}\phi_{s_{\mu,\rho}}(z)+(1+K_{1}e^{-\delta(t-T^{\ast})})\sigma\delta
e^{-\delta(t-T^{\ast})}\phi_{s_{\mu,\rho}}'(z)+F(\overline{u})\\
&\ \ +(1+K_{1}e^{-\delta(t-T^{\ast})})\big(s_{\mu,\rho}\phi_{s_{\mu,\rho}}'(z)-d_1 \phi_{s_{\mu,\rho}}''(z)+a\phi_{s_{\mu,\rho}}(z)
-b\psi_{s_{\mu,\rho}}(z)\big)\\
=&-\delta K_{1}e^{-\delta(t-T^{\ast})}\phi_{s_{\mu,\rho}}(z)+(1+K_{1}e^{-\delta(t-T^{\ast})})\sigma\delta
e^{-\delta(t-T^{\ast})}\phi_{s_{\mu,\rho}}'(z)\\
&\ \ +F(\overline{u})-(1+K_{1}e^{-\delta(t-T^{\ast})})F(\phi_{s_{\mu,\rho}}(z)).
\end{split}
\end{equation}
Note that $\overline{u}=(1+K_{1}e^{-\delta(t-T^{\ast})})\phi_{s_{\mu,\rho}}(z)>\phi_{s_{\mu,\rho}}(z)>0$  for $z\in(0,\overline{h}(t)-g(t))$. Hence
\begin{equation}\label{*22}
F(\overline{u})-(1+K_{1}e^{-\delta(t-T^{\ast})})F(\phi_{s_{\mu,\rho}}(z))
=\overline{u}\left[\frac{F(\overline{u})}{\overline{u}}-\frac{F(\phi_{s_{\mu,\rho}}(z))}{\phi_{s_{\mu,\rho}}(z)}\right]
\text{ for }z\in(0,\overline{h}(t)-g(t)).
\end{equation}

In the case when  $z\in(0,1]$, we define
$$\eta_1:=\mathop{\min}_{z \in [0,1]} {\phi_{s_{\mu,\rho}}'(z)}.$$
Then $\mathop{\inf}_{z \in (0,1]} {\phi_{s_{\mu,\rho}}'(z)}\ge \eta_1>0$. It follows from \eqref{ad ee}, \eqref{*22}, and  $\mathrm{(A2)}$ that
\begin{equation*}
\begin{split}
&~\overline{u}_{t}-d_1 \overline{u}_{xx}+a\overline{u}-b\overline{v}+F(\overline{u})\\
\geq&-\delta K_{1}e^{-\delta(t-T^{\ast})}\phi_{s_{\mu,\rho}}(z)+(1+K_{1}e^{-\delta(t-T^{\ast})})\sigma\delta
e^{-\delta(t-T^{\ast})}\phi_{s_{\mu,\rho}}'(z)\\
\geq&~e^{-\delta(t-T^{\ast})}\delta[-K_{1}u^*+\sigma\phi_{s_{\mu,\rho}}'(z)]\\
\geq& e^{-\delta(t-T^{\ast})}\delta(-K_{1}u^*+\sigma \eta_1)\ \text{ for } z\in(0,1].
\end{split}
\end{equation*}

In the case when $z>1$, we have $\overline{u}=(1+K_{1}e^{-\delta(t-T^{\ast})})\phi_{s_{\mu,\rho}}(z)>\phi_{s_{\mu,\rho}}(z)>\phi_{s_{\mu,\rho}}(1)$. In view of \eqref{*22}, it follows from $\mathrm{(A2)}$ that there exists a constant $\varrho_{1}>0$ small enough such that
\begin{equation}\label{*222}
F(\overline{u})-(1+K_{1}e^{-\delta(t-T^{\ast})})F(\phi_{s_{\mu,\rho}}(z)) =\overline{u}\left[\frac{F(\overline{u})}{\overline{u}}-\frac{F(\phi_{s_{\mu,\rho}}(z))}{\phi_{s_{\mu,\rho}}(z)}\right]
\ge\overline{u}\varrho_{1}\ \text{ for } z>1.
\end{equation}
In light of $\overline{u}\le(1+K_1)u^*$, it follows from \eqref{ad ee},   \eqref{*222}, and the strict monotonicity of $\phi_{s_{\mu,\rho}}$ that
\begin{equation*}
\begin{split}
&~\overline{u}_{t}-d_1 \overline{u}_{xx}+a\overline{u}-b\overline{v}+F(\overline{u})\\
\geq&-\delta K_{1}e^{-\delta(t-T^{\ast})}\phi_{s_{\mu,\rho}}(z)+F(\overline{u})-(1+K_{1}e^{-\delta(t-T^{\ast})})F(\phi_{s_{\mu,\rho}})\\
=&-\delta(\overline{u}-\phi_{s_{\mu,\rho}}(z))+\overline{u}\left[\frac{F(\overline{u})}{\overline{u}}
-\frac{F(\phi_{s_{\mu,\rho}}(z))}{\phi_{s_{\mu,\rho}}(z)}\right]\\
=&~\overline{u}\left[\frac{F(\overline{u})}{\overline{u}}-
\frac{F(\phi_{s_{\mu,\rho}}(z))}{\phi_{s_{\mu,\rho}}(z)}-\delta
+\delta\frac{\phi_{s_{\mu,\rho}}(z)}{\overline{u}}\right]\\
\geq&~
\overline{u}\left[\varrho_{1}-\delta+\delta\frac{\phi_{s_{\mu,\rho}}(1)}{(1+K_1)u^*}\right]\ \text{ for } z>1.
\end{split}
\end{equation*}
In conclusion, \eqref{3.2} holds provided that
\begin{equation}\label{3.10}
\sigma\geq \frac{K_{1}u^*}{\eta_1}\text{ and }\delta\le\frac{\varrho_{1}(1+K_1)u^*}{(1+K_1)u^*-\phi_{s_{\mu,\rho}}(1)}.
\end{equation}

Similarly, we define
$$\eta_2:=\mathop{\min}_{z \in [0,1]}{\psi_{s_{\mu,\rho}}'(z)}.$$
Then $\eta_2>0$ and $\overline{v}_{t}-d_2 \overline{v}_{xx}-c\overline{u}+d\overline{v}+G(\overline{v})\ge0$ for $z\in(0,1]$ provided that $\sigma \geq K_{1}v^*/\eta_2$. Moreover, in the case when $z>1$,  there exists a constant $\varrho_{2}>0$ such that $$\frac{G(\overline{v})}{\overline{v}}-\frac{G(\psi_{s_{\mu,\rho}})}{\psi_{s_{\mu,\rho}}}>\varrho_{2}\ \text{ for }z>1.$$
By parallel arguments, we deduce that \eqref{3.3} holds  provided that
\begin{equation}\label{3.11}
\sigma\geq \frac{K_{1}v^*}{\eta_2}\text{ and }\delta\le\frac{\varrho_{2}(1+K_1)v^*}{(1+K_1)v^*-\psi_{s_{\mu,\rho}}(1)}.
\end{equation}

Therefore, $(\overline{u},\overline{v},\overline{g},\overline{h})$ is an upper solution for system \eqref{1.2} when  \eqref{3.8}, \eqref{3.9}, \eqref{3.10}, and \eqref{3.11} hold simultaneously, which can be accomplished by  choosing suitable $K_1$, $\delta$, and $\sigma$. For example,
 one can choose
\begin{align*}
&K_1=\frac{C_0}{\min\{\phi_{s_{\mu,\rho}}(X_0),\psi_{s_{\mu,\rho}}(X_0)\}}-1,\\
&\delta=\min\left\{\frac{\varrho_{1}(1+K_1)u^*}{(1+K_1)u^*-\phi_{s_{\mu,\rho}}(1)},
\frac{\varrho_{2}(1+K_1)v^*}{(1+K_1)v^*-\psi_{s_{\mu,\rho}}(1)}\right\},\\
&\sigma=K_1\cdot\max\left\{\frac{s_{\mu,\rho}}{\delta},\frac{u^*}{\eta_1},\frac{v^*}{\eta_2}\right\}.
\end{align*}
It then follows from Lemma  \ref{lem2.1} that the desired result holds.
The proof is complete.
\end{proof}

We next construct a lower solution $(\underline{u},\underline{v},\underline{g},\underline{h})$ of system \eqref{1.2}. Define
\begin{equation}\label{lower sol}
\begin{split}
&\underline{h}(t)=s_{\mu,\rho}(t-T_{\ast})+L_0-(1-e^{-\delta(t-T_{\ast})}),\quad\underline{g}(t)
=-\underline{h}(t),\\
&\underline{u}(x,t)=(1-\tilde{\epsilon}e^{-\delta(t-T_{\ast})})
\left[\phi_{s_{\mu,\rho}}(\underline{h}(t)-x)+\phi_{s_{\mu,\rho}}(\underline{h}(t)+x)
-\phi_{s_{\mu,\rho}}(2\underline{h}(t))\right],\\
&\underline{v}(x,t)=(1-\tilde{\epsilon}e^{-\delta(t-T_{\ast})})
\left[\psi_{s_{\mu,\rho}}(\underline{h}(t)-x)+\psi_{s_{\mu,\rho}}(\underline{h}(t)+x)
-\psi_{s_{\mu,\rho}}(2\underline{h}(t))\right],
\end{split}
\end{equation}
where $(\phi_{s_{\mu,\rho}},\psi_{s_{\mu,\rho}})$ is the strictly increasing solution of system \eqref{1.5} with $s=s_{\mu,\rho}$, and $\delta$, $\tilde{\epsilon}\in(0,1)$, $T_{\ast}$ and $L_0>1$  are positive constants to be determined.

To verify that $(\underline{u},\underline{v},\underline{g},\underline{h})$ is a lower solution of system \eqref{1.2}, the following result is needed.

\begin{lem}\label{lem3.3}
{\rm(\cite[Lemma 2.9]{DM2022})}
Assume that $\mathbf{F}=(f_i) \in C^{2}(\mathbb{R}^{m}, \mathbb{R}^{m})$, $\mathbf{u^*}=(u_1^*,\cdots,u_m^*)\succ\mathbf{0}$ and
\begin{equation*}
\mathbf{F}(\mathbf{u^*})=\mathbf{0},~\mathbf{u^*}[\nabla \mathbf{F}(\mathbf{u^*})]^\intercal\prec\mathbf{0}.
\end{equation*}
Then there exists $\delta_{0}>0$ small enough such that for $0<\epsilon\ll 1$ and $\mathbf{u},\mathbf{v}\in[(1-\delta_0)\mathbf{u^*},\mathbf{u^*}]$ satisfying
$$(u^*_{i}-u_{i})(u^*_{j}-v_{j})\leq C\delta_{0}\epsilon\text{ for some }C>0\text{ and all }i,j \in \{1,\cdots,m\},$$
we have
$$(1-\epsilon)[\mathbf{F}(\mathbf{u})+\mathbf{F}(\mathbf{v})]-\mathbf{F}((1-\epsilon)(\mathbf{u}+\mathbf{v}-\mathbf{u^*}))\preceq\frac{\epsilon}{2}
\mathbf{u^*}[\nabla \mathbf{F}(\mathbf{u^*})]^\intercal,$$
where $\mathbf{u}=(u_1,\cdots,u_m)$ and $\mathbf{v}=(v_1,\cdots,v_m)$.
\end{lem}

We will use this lemma with $\mathbf{u^*}=(u^*,v^*)$ and $\mathbf{F}(u,v)=f_i(u,v)$ for $i=1,2$, where
$$f_{1}(u,v):=-au+bv-F(u)\text{ and }f_{2}(u,v):=cu-dv-G(v).$$
\begin{lem}\label{lem3.4}
Assume that $\mathrm{(H)}$ holds and let $(u,v,g,h)$ be the unique solution of system \eqref{1.2}. If spreading occurs,  then there exist suitable positive constants  $\delta$, $\tilde{\epsilon}\in(0,1)$, $T_{\ast}$, and $L_0>1$ such that
$$h(t)\ge\underline{h}(t),~(u(x,t),v(x,t))\succeq(\underline{u}(x,t),\underline{v}(x,t))\text{ for }[\underline{g}(t),\underline{h}(t)],~t>T_*.$$
\end{lem}
\begin{proof}
We claim that  $(\underline{u},\underline{v},\underline{g},\underline{h})$ is a lower solution for system \eqref{1.2} by taking  appropriate parameters  $\delta$, $\tilde{\epsilon}$, $T_{\ast}$, and $L_0$, that is,
\begin{align}
&\underline{u}_{t}\leq d_1 \underline{u}_{xx}+f_{1}(\underline{u},\underline{v}), &&\underline{g}(t)<x<\underline{h}(t),~t>T_{\ast},\label{3.14}\\
&\underline{v}_{t}\leq d_2 \underline{v}_{xx}+f_{2}(\underline{u},\underline{v}),
&&\underline{g}(t)<x<\underline{h}(t),~t>T_{\ast},\label{3.15}\\
&\underline{u}(x,t)=\underline{v}(x,t)=0,&&x=\underline{g}(t)\text{ or }x=\underline{h}(t),~t>T_{\ast},\label{3.17}\\
&\underline{g}(T_{\ast})\geq g(T_{\ast}),~\underline{g}'(t)\geq -\mu [\underline{v}_{x}(\underline{g}(t),t)+\rho\underline{u}_{x}(\underline{g}(t),t)],&&t>T_{\ast},\label{3.18.2}\\
&\underline{h}(T_{\ast})\leq h(T_{\ast}),~\underline{h}'(t)\leq -\mu [\underline{v}_{x}(\underline{h}(t),t)+\rho\underline{u}_{x}(\underline{h}(t),t)],&&t>T_{\ast},\label{3.18}\\
&\underline{u}(x,T_{\ast})\leq {u}(x,T_{\ast}),~\underline{v}(x,T_{\ast})\leq {v}(x,T_{\ast}),&&\underline{g}(T_{\ast})\leq x\leq\underline{h}(T_{\ast}).\label{3.19}
\end{align}

First, we have
$(\underline{u}(\underline{g}(t),t),\underline{v}(\underline{g}(t),t))=
 (\underline{u}(\underline{h}(t),t),\underline{v}(\underline{h}(t),t))=(0,0)$ for $t> T^*$. This proves \eqref{3.17}.

We next show \eqref{3.18.2} and \eqref{3.18}. Since spreading occurs, there exists $T_*>0$ large enough such that
\begin{equation*}
[\underline{g}(T_*),\underline{h}(T_*)]=[-L_0,L_0]\subset[g(T_*),h(T_*)].
\end{equation*}
Moreover, a direct calculation yields $\underline{h}'(t)=s_{\mu,\rho}-\delta e^{-\delta(t-T_{\ast})}$ for $t>T_*$. It follows from \eqref{s mu} that for $t>T_*$,
\begin{equation*}
\begin{split}
&-\mu[\underline{v}_{x}(\underline{h}(t),t)+\rho\underline{u}_{x}(\underline{h}(t),t)]\\
=&~\mu(1-\tilde{\epsilon}e^{-\delta(t-T_{\ast})})
\big[\psi_{s_{\mu,\rho}}'(0)+\rho\phi_{s_{\mu,\rho}}'(0)-\psi_{s_{\mu,\rho}}'(2\underline{h}(t))-\rho\phi_{s_{\mu,\rho}}'(2\underline{h}(t))\big]\\
=&~s_{\mu,\rho}(1-\tilde{\epsilon}e^{-\delta(t-T_{\ast})})-\mu(1-\tilde{\epsilon}e^{-\delta(t-T_{\ast})})
\big[\psi_{s_{\mu,\rho}}'(2\underline{h}(t))+\rho\phi_{s_{\mu,\rho}}'(2\underline{h}(t))\big].
\end{split}
\end{equation*}
According to  Lemma \ref{lem2.4}, there exist constants $C, \hat{\mu}_{1}>0$ such that for $L_0\gg 1$ and $t>T^*$,
$$\phi_{s_{\mu,\rho}}'(2\underline{h}(t))\le Ce^{-2\hat{\mu}_1\underline{h}(t)}\text{ and }
\psi_{s_{\mu,\rho}}'(2\underline{h}(t))\le  Ce^{-2\hat{\mu}_1\underline{h}(t)}.$$
Then  one can find a large constant $C_1>0$ such that for $t>T^*$,
\begin{equation*}
\begin{split}
&-\mu[\underline{v}_{x}(\underline{h}(t),t)+\rho\underline{u}_{x}(\underline{h}(t),t)]\\
\geq&~s_{\mu,\rho}-s_{\mu,\rho}\tilde{\epsilon}e^{-\delta(t-T_{\ast})}-C_{1}e^{-2\hat{\mu}_{1}\underline{h}(t)}\\
\geq&~s_{\mu,\rho}-s_{\mu,\rho}\tilde{\epsilon}e^{-\delta(t-T_{\ast})}-C_{1}e^{-2\hat{\mu}_{1}s_{\mu,\rho}(t-T_{\ast})}e^{-2\hat{\mu}_{1}(L_0-1)}\\
\geq&~s_{\mu,\rho}-e^{-\delta(t-T_{\ast})}\big(s_{\mu,\rho}\tilde{\epsilon}+C_{1}e^{-2\hat{\mu}_{1}(L_0-1)}\big),
\end{split}
\end{equation*}
where the last inequality holds when
\begin{equation*}\label{delta}
\delta\in(0,2\hat{\mu}_{1}s_{\mu,\rho}).
\end{equation*}
 Thus, $\underline{h}'(t)\le -\mu[\underline{v}_{x}(\underline{h}(t),t)+\rho\underline{u}_{x}(\underline{h}(t),t)]$ for $t>T_*$ provided that \begin{equation}\label{3.20}
s_{\mu,\rho}\tilde{\epsilon}+C_{1}e^{-2\hat{\mu}_{1}(L_0-1)}\leq\delta.
\end{equation}
Since $\underline{u}(x,t)$ and $\underline{v}(x,t)$ are even in $x$, and $\underline{g}(t)=-\underline{h}(t)$, we obtain  $\underline{g}'(t)\ge -\mu[\underline{v}_{x}(\underline{g}(t),t)+\rho\underline{u}_{x}(\underline{g}(t),t)]$ when
 $\underline{h}'(t)\le -\mu[\underline{v}_{x}(\underline{h}(t),t)+\rho\underline{u}_{x}(\underline{h}(t),t)]$ holds. Therefore, both \eqref{3.18.2} and \eqref{3.18} are valid when \eqref{3.20} holds.

 Since spreading occurs, one can enlarge $T_{\ast}$ if necessary  and choose $\tilde{\epsilon}\in(0,1)$ small enough such that
 $$((1-\tilde{\epsilon})u^*,(1-\tilde{\epsilon})v^*)
\preceq(u(x,T_{\ast}),v(x,T_{\ast}))\text{ for }x\in[\underline{g}(T_{\ast}), \underline{h}(T_{\ast})]=[-L_0,L_0].$$
By the monotonicity of  $\phi_{s_{\mu,\rho}}$ and $\psi_{s_{\mu,\rho}}$, one can choose $\tilde{\epsilon}\in(0,1)$ small to deduce that
\begin{equation*}
\begin{split}
&\underline{u}(x,T_{\ast})\le (1-\tilde{\epsilon})\phi_{s_{\mu,\rho}}(L_0-x)\le (1-\tilde{\epsilon})u^*\text{ for }x\in[-L_0,L_0],\\
&\underline{v}(x,T_{\ast})\le (1-\tilde{\epsilon})\psi_{s_{\mu,\rho}}(L_0-x)\le (1-\tilde{\epsilon})v^*\text{ for }x\in[-L_0,L_0],
\end{split}
\end{equation*}
which leads to $(\underline{u}(x,T_{\ast}),\underline{v}(x,T_{\ast}))\preceq(u(x,T_{\ast}),v(x,T_{\ast}))$ for $x\in[-L_0,L_0]$. This  means that \eqref{3.19} holds.

It remains to check \eqref{3.14} and \eqref{3.15}. For convenience, we set
$$\underline{z}:=\underline{h}(t)-x,~\overline{z}:=\underline{h}(t)+x,~\epsilon(t)=\tilde{\epsilon}e^{-\delta(t-T_{\ast})}.$$
Then
$$\underline{z},\ \overline{z}\in(0,2\underline{h}(t))\text{ for }x\in(-\underline{h}(t),\underline{h}(t)),\ t>T_*$$
 and
\begin{equation}\label{uv}
\begin{split}
&\underline{u}(x,t)=(1-\epsilon(t))\big[\phi_{s_{\mu,\rho}}(\underline{z})+\phi_{s_{\mu,\rho}}(\overline{z})-\phi_{s_{\mu,\rho}}(2\underline{h}(t))\big],\\
&\underline{v}(x,t)=(1-\epsilon(t))\big[\psi_{s_{\mu,\rho}}(\underline{z})+\psi_{s_{\mu,\rho}}(\overline{z})-\psi_{s_{\mu,\rho}}(2\underline{h}(t))\big].
\end{split}
\end{equation}
In view of $\epsilon'(t)=-\delta\tilde{\epsilon}e^{-\delta(t-T_{\ast})}=-\delta\epsilon(t)$ and $\underline{h}'(t)=s_{\mu,\rho}-\delta e^{-\delta(t-T_{\ast})}$, we have
\begin{equation*}
\begin{split}
&~\underline{u}_{t}
=(1-\epsilon(t))(s_{\mu,\rho}-\delta e^{-\delta(t-T_{\ast})})\big[\phi_{s_{\mu,\rho}}'(\underline{z})
+\phi_{s_{\mu,\rho}}'(\overline{z})-2\phi_{s_{\mu,\rho}}'(2\underline{h}(t))\big]\\
&\quad\quad\ +\delta\epsilon(t)\big[\phi_{s_{\mu,\rho}}(\underline{z})
+\phi_{s_{\mu,\rho}}(\overline{z})-\phi_{s_{\mu,\rho}}(2\underline{h}(t))\big],\\
&~\underline{v}_{t}
=(1-\epsilon(t))(s_{\mu,\rho}-\delta e^{-\delta(t-T_{\ast})})\big[\psi_{s_{\mu,\rho}}'(\underline{z})
+\psi_{s_{\mu,\rho}}'(\overline{z})-2\psi_{s_{\mu,\rho}}'(2\underline{h}(t))\big]\\
&\quad\quad\ +\delta\epsilon(t)\big[\psi_{s_{\mu,\rho}}(\underline{z})
+\psi_{s_{\mu,\rho}}(\overline{z})-\psi_{s_{\mu,\rho}}(2\underline{h}(t))\big],\\
&~\underline{u}_{xx}=(1-\epsilon(t))\big[\phi_{s_{\mu,\rho}}''(\underline{z})+
\phi_{s_{\mu,\rho}}''(\overline{z})\big],\\
&~\underline{v}_{xx}=(1-\epsilon(t))\big[\psi_{s_{\mu,\rho}}''(\underline{z})+
\psi_{s_{\mu,\rho}}''(\overline{z})\big].
\end{split}
\end{equation*}
Due to $\phi_{s_{\mu,\rho}}<u^*$ and $\phi_{s_{\mu,\rho}}'>0$,  it follows from the equation of $\phi_{s_{\mu,\rho}}$ (see \eqref{1.5})   that for $x\in(-\underline{h}(t),\underline{h}(t))$ and $t>T_*$,
\begin{equation}\label{eq*1}
\begin{split}
&\qquad\underline{u}_{t}-d_1 \underline{u}_{xx}-f_{1}(\underline{u},\underline{v})\\
&=(1-\epsilon(t))(s_{\mu,\rho}-\delta e^{-\delta(t-T_{\ast})})\big[\phi_{s_{\mu,\rho}}'(\underline{z})
+\phi_{s_{\mu,\rho}}'(\overline{z})-2\phi_{s_{\mu,\rho}}'(2\underline{h}(t))\big]-f_{1}(\underline{u},\underline{v})\\
&\quad\ +\delta\epsilon(t)\big[\phi_{s_{\mu,\rho}}(\underline{z})
+\phi_{s_{\mu,\rho}}(\overline{z})-\phi_{s_{\mu,\rho}}(2\underline{h}(t))\big]-d_1(1-\epsilon(t))\big[\phi_{s_{\mu,\rho}}''(\underline{z})+
\phi_{s_{\mu,\rho}}''(\overline{z})\big]\\
&=(1-\epsilon(t))\big[f_{1}(\phi_{s_{\mu,\rho}}(\underline{z}),\psi_{s_{\mu,\rho}}(\underline{z}))
+f_{1}(\phi_{s_{\mu,\rho}}(\overline{z}),\psi_{s_{\mu,\rho}}(\overline{z}))\big]
-f_{1}(\underline{u},\underline{v})\\
&\quad\ +\delta\epsilon(t)\big[\phi_{s_{\mu,\rho}}(\underline{z})
+\phi_{s_{\mu,\rho}}(\overline{z})-\phi_{s_{\mu,\rho}}(2\underline{h}(t))\big]-(1-\epsilon(t))\delta e^{-\delta(t-T_{\ast})}\big[\phi_{s_{\mu,\rho}}'(\underline{z})+\phi_{s_{\mu,\rho}}'(\overline{z})\big]\\
&\quad\ -
2(1-\epsilon(t))(s_{\mu,\rho}-\delta e^{-\delta(t-T_{\ast})})\phi_{s_{\mu,\rho}}'(2\underline{h}(t))\\
&\le(1-\epsilon(t))[f_{1}(\phi_{s_{\mu,\rho}}(\underline{z}),\psi_{s_{\mu,\rho}}(\underline{z}))
+f_{1}(\phi_{s_{\mu,\rho}}(\overline{z}),\psi_{s_{\mu,\rho}}(\overline{z}))]
-f_{1}(\underline{u},\underline{v})\\
&\quad\ +2\delta \epsilon (t)u^*-(1-\epsilon(t))\delta e^{-\delta(t-T_{\ast})}\big[\phi_{s_{\mu,\rho}}'(\underline{z})+\phi_{s_{\mu,\rho}}'(\overline{z})\big]\\
&=B_1(x,t)+B_2(x,t)+B_3(x,t)+2\delta \epsilon (t)u^*
\end{split}
\end{equation}
when $\delta\in(0,s_{\mu,\rho})$, where
\begin{equation*}
\begin{split}
&B_1(x,t):=(1-\epsilon(t))\big[f_{1}(\phi_{s_{\mu,\rho}}(\underline{z}),\psi_{s_{\mu,\rho}}(\underline{z}))
         +f_{1}(\phi_{s_{\mu,\rho}}(\overline{z}),\psi_{s_{\mu,\rho}}(\overline{z}))\big]\\
&\qquad\qquad\ \ -f_{1}\big((1-\epsilon(t))
(\phi_{s_{\mu,\rho}}(\underline{z})+\phi_{s_{\mu,\rho}}(\overline{z})-u^*),(1-\epsilon(t))
(\psi_{s_{\mu,\rho}}(\underline{z})+\psi_{s_{\mu,\rho}}(\overline{z})-v^*)\big),\\
&B_2(x,t):=f_{1}\big((1-\epsilon(t))
(\phi_{s_{\mu,\rho}}(\underline{z})+\phi_{s_{\mu,\rho}}(\overline{z})-u^*),(1-\epsilon(t))
(\psi_{s_{\mu,\rho}}(\underline{z})+\psi_{s_{\mu,\rho}}(\overline{z})-v^*)\big)\\
&\qquad\qquad\ \ -f_{1}(\underline{u},\underline{v}),\\
&B_3(x,t):=-(1-\epsilon(t))\delta e^{-\delta(t-T_{\ast})}\big[\phi_{s_{\mu,\rho}}'(\underline{z})+\phi_{s_{\mu,\rho}}'(\overline{z})\big].
\end{split}
\end{equation*}

In light of  \eqref{uv}, we have
\begin{align*}\label{eq*01}
&(1-\epsilon(t))
(\phi_{s_{\mu,\rho}}(\underline{z})+\phi_{s_{\mu,\rho}}(\overline{z})-u^*)=
\underline{u}+(1-\epsilon(t))(\phi_{s_{\mu,\rho}}(2\underline{h}(t))-u^*),\\
&(1-\epsilon(t))
(\psi_{s_{\mu,\rho}}(\underline{z})+\psi_{s_{\mu,\rho}}(\overline{z})-v^*)=
\underline{v}+(1-\epsilon(t))(\psi_{s_{\mu,\rho}}(2\underline{h}(t))-v^*).
\end{align*}
Since $f_1$ is Lipschitz continuous, there exist  constants $L_1, L_2>0$  such that for $x\in(-\underline{h}(t),\underline{h}(t))$ and $t>T^*$,
\begin{equation}\label{eq*2}
\begin{split}
&\quad\ B_2(x,t)\\
&=f_1\big(\underline{u}+(1-\epsilon(t))(\phi_{s_{\mu,\rho}}(2\underline{h}(t))-u^*),
\underline{v}+(1-\epsilon(t))(\psi_{s_{\mu,\rho}}(2\underline{h}(t))-v^*)\big) -f_1(\underline{u},\underline{v})\\
&\le L_1(1-\epsilon(t))\big(u^*-\phi_{s_{\mu,\rho}}(2\underline{h}(t))\big)+L_2(1-\epsilon(t))
\big(v^*-\psi_{s_{\mu,\rho}}(2\underline{h}(t))\big)\\
&\le L_1\big(u^*-\phi_{s_{\mu,\rho}}(2\underline{h}(t))\big)+
L_2\big(v^*-\psi_{s_{\mu,\rho}}(2\underline{h}(t))\big).
\end{split}
\end{equation}
 It follows from Lemma \ref{lem2.4} that there exist constants $p, q, \hat{\mu}_1>0$ such that  for   $L_0\gg1$ and $t>T_*$,
\begin{equation*}
\phi_{s_{\mu,\rho}}(2\underline{h}(t))\ge u^*-pe^{-2\hat{\mu}_1\underline{h}(t)}\text{ and }
\psi_{s_{\mu,\rho}}(2\underline{h}(t))\ge v^*-qe^{-2\hat{\mu}_1\underline{h}(t)}.
\end{equation*}
Combining with \eqref{eq*2}, we conclude that
 \begin{equation}\label{eq*000}
B_2(x,t)\le C_2e^{-2\hat{\mu}_1\underline{h}(t)}\text{ for }x\in(-\underline{h}(t),\underline{h}(t)),~t>T_*,
\end{equation}
 where $C_2:=L_1p+L_2q>0$.

Therefore, it follows from \eqref{eq*1} and \eqref{eq*000} that for $x\in(-\underline{h}(t),\underline{h}(t))$ and $t>T_*$,
\begin{equation}\label{eq*3}
\underline{u}_{t}-d_1 \underline{u}_{xx}-f_{1}(\underline{u},\underline{v})\le
2\delta \epsilon (t)u^*+C_2e^{-2\hat{\mu}_1\underline{h}(t)}+B_1(x,t)+B_3(x,t).
\end{equation}
In subsequent discussions, we divide the interval $(-\underline{h}(t),\underline{h}(t))$ into the following three subintervals:
\begin{equation*}
\begin{split}
&I_{1}(t):=(-\underline{h}(t),-\underline{h}(t)+K_{0}],\\
&I_{2}(t):=(-\underline{h}(t)+K_{0},\underline{h}(t)-K_{0}),\\
&I_{3}(t):=[\underline{h}(t)-K_{0},\underline{h}(t)),
\end{split}
\end{equation*}
where $K_0>0$ is a given constant.

\textbf{Claim 1:} There exist constants $L_0\gg 1$, $0<\tilde{\epsilon}\ll 1$, and $\delta\in(0,\min\{s_{\mu,\rho},2\hat{\mu}_1s_{\mu,\rho})\}$ such that
\begin{equation*}\label{eq*444}
\underline{u}_{t}-d_1 \underline{u}_{xx}-f_{1}(\underline{u},\underline{v})\le 0\text{ for }x\in I_1(t)\cup I_3(t),~t>T_*.
\end{equation*}

To prove Claim 1, we first show that there exists a constant $\eta_0>0$ such that  for $x\in I_1(t)\cup I_3(t)$ and $t>T_*$,
\begin{equation}\label{eq*4}
B_1(x,t)=O(\epsilon(t))+O(e^{-2\hat{\mu}_1\underline{h}(t)})\text{ and }
B_3(x,t)\le -(1-\tilde{\epsilon})\delta\eta_0 e^{-\delta(t-T_*)}.
\end{equation}
Since both $B_1(x,t)$ and $B_3(x,t)$ are symmetric in $x$, and  $x\in I_1(t)$ if and only if $-x\in I_3(t)$,  it suffices to prove \eqref{eq*4} when $x\in I_3(t)$ and $t>T_*$. Fix $x\in I_3(t)=[\underline{h}(t)-K_0,\underline{h}(t))$ and $t>T_*$. Then
$$\underline{z}=\underline{h}(t)-x\in (0,K_0]~\text{ and }~\overline{z}=\underline{h}(t)+x\in[2\underline{h}(t)-K_0,2\underline{h}(t)).$$
Define
$$\eta_0:=\min\limits_{z\in[0,K_0]}\phi_{s_{\mu,\rho}}'(z).$$
Then  $\inf_{z\in(0,K_0]}\phi_{s_{\mu,\rho}}'(z)\ge\eta_0>0$, which yields
$$\phi_{s_{\mu,\rho}}'(\underline{z})+\phi_{s_{\mu,\rho}}'(\overline{z})\ge\phi_{s_{\mu,\rho}}'(\underline{z})
\ge\eta_0\text{ for }x\in I_3(t),~t>T_*.$$
Combining with the monotonicity of $\phi_{s_{\mu,\rho}}$, we obtain
\begin{equation}\label{eq*5}
B_3(x,t)\le -(1-\tilde{\epsilon})\delta\eta_0e^{-\delta(t-T_*)}\text{ for }x\in I_3(t),\ t>T_*.
\end{equation}
 It follows from  Lemma  \ref{lem2.4} that for $L_0\gg 1$,
\begin{equation*}\label{aequation 01}
\phi_{s_{\mu,\rho}}(\overline{z})=u^*+O(e^{-2\hat{\mu}_1\underline{h}(t)})\text{ and }\psi_{s_{\mu,\rho}}(\overline{z})=v^*+O(e^{-2\hat{\mu}_1\underline{h}(t)}),
\end{equation*}
which implies that
$$B_1(x,t)=O(\epsilon(t))+O(e^{-2\hat{\mu}_1\underline{h}(t)}).$$
Combining with \eqref{eq*5}, we get \eqref{eq*4}.
 Hence  it follows from \eqref{eq*3} and \eqref{eq*4} that there exist constants $C_3, C_4>0$ such that
for $x\in I_1(t)\cup I_3(t)$ and $t>T_*$,
\begin{equation*}
\begin{split}
~&\qquad\underline{u}_{t}-d_1 \underline{u}_{xx}-f_{1}(\underline{u},\underline{v})\\
~&\le 2\delta \epsilon (t)u^*+C_2e^{-2\hat{\mu}_1\underline{h}(t)}+O(\epsilon(t))+O(e^{-2\hat{\mu}_1\underline{h}(t)})
-(1-\tilde{\epsilon})\delta\eta_0 e^{-\delta(t-T_*)}\\
~&\le C_3\tilde{\epsilon}e^{-\delta(t-T_*)}+C_4e^{-2\hat{\mu}_1\underline{h}(t)}-\delta\eta_0 e^{-\delta(t-T_*)}\\
~&\le C_3\tilde{\epsilon}e^{-\delta(t-T_*)}+C_4e^{-2\hat{\mu}_1s_{\mu,\rho}(t-T_*)}e^{-2\hat{\mu}_1(L_0-1)}-\delta\eta_0 e^{-\delta(t-T_*)}\\
~&\le  \big(C_3\tilde{\epsilon}+C_4e^{-2\hat{\mu}_1(L_0-1)}-\delta\eta_0\big)e^{-\delta(t-T_*)}\le0
\end{split}
\end{equation*}
provided that
\begin{equation}\label{condition 1}
\delta\in(0,\min\{s_{\mu,\rho},2\hat{\mu}_1s_{\mu,\rho})\})\text{ and }C_3\tilde{\epsilon}+C_4e^{-2\hat{\mu}_1(L_0-1)}\le\delta\eta_0,
\end{equation}
which can be guaranteed by choosing $L_0\gg1$  and $0<\tilde{\epsilon}\ll1$. This proves Claim 1.

\textbf{Claim 2:} There exist constants  $L_0\gg 1$, $0<\tilde{\epsilon}\ll 1$, and $\delta\in(0,\min\{s_{\mu,\rho},\hat{\mu}_1s_{\mu,\rho})\})$ such that
$$\underline{u}_{t}-d_1 \underline{u}_{xx}-f_{1}(\underline{u},\underline{v})\le 0\text{ for }x\in I_2(t),~t>T_*.$$

Due to $\phi_{s_{\mu,\rho}}'>0$ and $\epsilon(t)<1$ for $t>T_*$, we have 
$$B_3(x,t):=-(1-\epsilon(t))\delta e^{-\delta(t-T_{\ast})}\big[\phi_{s_{\mu,\rho}}'(\underline{z})+\phi_{s_{\mu,\rho}}'(\overline{z})\big]
\le0\text{ for }x\in I_2(t),\ t>T_*.$$
Then  \eqref{eq*3} yields
\begin{equation}\label{eq*6}
\underline{u}_{t}-d_1 \underline{u}_{xx}-f_{1}(\underline{u},\underline{v})\le
2\delta \epsilon (t)u^*+C_2e^{-2\hat{\mu}_1\underline{h}(t)}+B_1(x,t)\text{ for }x\in I_2(t),\ t>T_*.
\end{equation}
We next apply Lemma \ref{lem3.3} with
$$\mathbf{F}=f_1, \mathbf{u}^*=(u^*,v^*), \mathbf{u}=(\phi_{s_{\mu,\rho}}(\underline{z}),\psi_{s_{\mu,\rho}}(\underline{z})),
\text{ and }\mathbf{v}=(\phi_{s_{\mu,\rho}}(\overline{z}),\psi_{s_{\mu,\rho}}(\overline{z}))$$
to show that there exists $\sigma_0>0$ such that
\begin{equation}\label{eq*7}
B_1(x,t)\le-\sigma_0\epsilon(t)\text{ for }x\in I_2(t),\ t>T_*.
\end{equation}
Obviously, we have $f_1(u^*,v^*)=-au^*+bv^*-F(u^*)=0$. Moreover, it follows from $\mathrm{(A2)}$ that
$$\left(\frac{F(\zeta)}{\zeta}\right)'=\frac{F'(\zeta)\zeta-F(\zeta)}{\zeta^2}>0\text{ for }\zeta>0,$$
which implies that
$$(u^*,v^*)[\nabla f_1(u^*,v^*)]^\intercal=(u^*,v^*)(-a-F'(u^*),b)^\intercal=F(u^*)-F'(u^*)u^*<0.$$

For $x\in I_2(t)=(-\underline{h}(t)+K_0,\underline{h}(t)-K_0)$ and $t>T_*$, we have
\begin{equation}\label{eq*8}
\underline{z},\ \overline{z}\in(K_0,2\underline{h}(t)-K_0).
\end{equation}
Thus, by virtue of $(\phi_{s_{\mu,\rho}}(+\infty),\psi_{s_{\mu,\rho}}(+\infty))=(u^*,v^*)$, we can choose $K_0>0$ large enough to find a small constant $\delta_0>0$ such that
\begin{equation}\label{eq*9}
u^*-\phi_{s_{\mu,\rho}}(K_0)\in(0, \delta_0],\ v^*-\psi_{s_{\mu,\rho}}(K_0)\in(0, \delta_0]
\end{equation}
and
\begin{equation*}
\phi_{s_{\mu,\rho}}(\underline{z}),\ \phi_{s_{\mu,\rho}}(\overline{z})\in[(1-\delta_0)u^*,u^*],\ \
\psi_{s_{\mu,\rho}}(\underline{z}),\ \psi_{s_{\mu,\rho}}(\overline{z})\in[(1-\delta_0)v^*,v^*]\text{ for }x\in I_2(t),\ t>T_*.
\end{equation*}

In the case when $x\in(-\underline{h}(t)+K_0,0]$, we have  $\underline{z}=\underline{h}(t)-x\ge\underline{h}(t)$. It follows from Lemma \ref{lem2.4} that there exists $C_5>0$ such that for $L_0\gg 1$,
$$u^*-\phi_{s_{\mu,\rho}}(\underline{h}(t)),\  v^*-\psi_{s_{\mu,\rho}}(\underline{h}(t))\le C_5 e^{-\hat{\mu}_1\underline{h}(t)}\text{ for }t>T_*.$$
 Combining with the monotonicity of $\phi_{s_{\mu,\rho}}$ and $\psi_{s_{\mu,\rho}}$, we have
\begin{equation}\label{eq*10}
\begin{split}
&0<u^*-\phi_{s_{\mu,\rho}}(\underline{z})\le u^*-\phi_{s_{\mu,\rho}}(\underline{h}(t))\le C_5 e^{-\hat{\mu}_1\underline{h}(t)}\text{ for }x\in I_2(t),\ t>T_*,\\
&0<v^*-\psi_{s_{\mu,\rho}}(\underline{z})\le v^*-\psi_{s_{\mu,\rho}}(\underline{h}(t))\le C_5 e^{-\hat{\mu}_1\underline{h}(t)}\text{ for }x\in I_2(t),\ t>T_*.
\end{split}
\end{equation}
In view of \eqref{eq*8} and \eqref{eq*9}, we use the monotonicity of $\phi_{s_{\mu,\rho}}$ and $\psi_{s_{\mu,\rho}}$ again to deduce that
\begin{equation*}
\begin{split}
&0<u^*-\phi_{s_{\mu,\rho}}(\overline{z})\le u^*-\phi_{s_{\mu,\rho}}(K_0)\le \delta_0\text{ for }x\in I_2(t),\ t>T_*,\\
&0<v^*-\psi_{s_{\mu,\rho}}(\overline{z})\le v^*-\psi_{s_{\mu,\rho}}(K_0)\le \delta_0\text{ for }x\in I_2(t),\ t>T_*.
\end{split}
\end{equation*}
Combining with \eqref{eq*10}, it follows that for $x\in I_2(t)$ and  $t>T_*$,
 \begin{equation}\label{eq*11}
\begin{split}
& (u^*-\phi_{s_{\mu,\rho}}(\underline{z}))(u^*-\phi_{s_{\mu,\rho}}(\overline{z}))\le \delta_0 C_5 e^{-\hat{\mu}_1\underline{h}(t)}\le  \delta_0 C_5 e^{-\hat{\mu}_1[s_{\mu,\rho}(t-T_*)+L_0-1]}\le\delta_0 C_5\epsilon(t),\\
& (v^*-\psi_{s_{\mu,\rho}}(\underline{z}))(v^*-\psi_{s_{\mu,\rho}}(\overline{z}))\le \delta_0 C_5 e^{-\hat{\mu}_1\underline{h}(t)}\le  \delta_0 C_5 e^{-\hat{\mu}_1[s_{\mu,\rho}(t-T_*)+L_0-1]}\le\delta_0 C_5\epsilon(t)
\end{split}
\end{equation}
provided that
 \begin{equation}\label{eq*12}
\delta\in(0,\hat{\mu}_1s_{\mu,\rho})\text{ and }e^{-\hat{\mu}_1(L_0-1)}\le\tilde{\epsilon}.
\end{equation}

In the case when $x\in[0,\underline{h}(t)-K_0)$, one can similarly obtain \eqref{eq*11} when \eqref{eq*12} holds. Therefore, it follows from Lemma \ref{lem3.3} that
$$ B_1(x,t)\le-\sigma_0\epsilon(t)\text{ for }x\in I_2(t),\ t>T_*$$
with $\sigma_0:=[F'(u^*)u^*-F(u^*)]/2>0$. This proves \eqref{eq*7}.

From \eqref{eq*6} and \eqref{eq*7}, we know that for $x \in I_{2}(t)$ and $t\geq T_{\ast}$,
\begin{equation*}
\begin{split}
&\quad\ \underline{u}_{t}-d_1 \underline{u}_{xx}-f_{1}(\underline{u},\underline{v})\\
&\le 2\delta \epsilon (t)u^*+C_2e^{-2\hat{\mu}_1\underline{h}(t)}-\sigma_0\epsilon(t)\\
&\leq2\delta\epsilon(t) u^*+C_{2}e^{-2\hat{\mu}_{1}s_{\mu,\rho}(t-T_{\ast})}e^{-2\hat{\mu}_{1}(L_0-1)}
-\sigma_{0}\epsilon(t)\\
&\leq[C_{2}e^{-2\hat{\mu}_{1}(L_0-1)}+(2\delta u^*-\sigma_{0})\tilde{\epsilon}]e^{-\delta(t-T_{\ast})}\le0
\end{split}
\end{equation*}
provided that
\begin{equation}\label{condition 2}
\delta\in(0,\min\{s_{\mu,\rho},2\hat{\mu}_1s_{\mu,\rho})\})\text{ and }
C_{2}e^{-2\hat{\mu}_{1}(L_0-1)}+(2\delta u^*-\sigma_{0})\tilde{\epsilon}\le0,
\end{equation}
which proves Claim 2.

Therefore, we conclude from  Claims 1 and 2 that
$$\underline{u}_{t}-d_1 \underline{u}_{xx}-f_{1}(\underline{u},\underline{v})\le 0\text{ for }x\in(-\underline{h}(t),\underline{h}(t))\text{ and }t>T_*$$
 when   \eqref{condition 1}, \eqref{eq*12}, and \eqref{condition 2} hold simultaneously, namely
\begin{equation}\label{Condition}
\begin{cases}
\delta\in(0,\min\{s_{\mu,\rho},\hat{\mu}_1s_{\mu,\rho})\}),\\
C_3\tilde{\epsilon}+C_4e^{-2\hat{\mu}_1(L_0-1)}\le\delta\eta_0,\\
e^{-\hat{\mu}_1(L_0-1)}\le\tilde{\epsilon},\\
C_{2}e^{-2\hat{\mu}_{1}(L_0-1)}+(2\delta u^*-\sigma_{0})\tilde{\epsilon}\le0.
\end{cases}
\end{equation}
 Since $e^{-\hat{\mu}_{1}(L_0-1)}\to+\infty$ exponentially as $L_0\to+\infty$,
one can fix $\delta\in(0,\min\{s_{\mu,\rho},\hat{\mu}_1s_{\mu,\rho}, \sigma_0/(2u^*)\})$,  and then  choose $0<\tilde{\epsilon}\ll 1$  and  $L_0\gg 1$  to guarantee \eqref{Condition}. Hence \eqref{3.14} holds. Similarly, one can choose suitable  $\delta, \tilde{\epsilon}$, and $L_0$   to prove  \eqref{3.15}.

In conclusion,  there exist suitable constants $\delta$, $\tilde{\epsilon}\in(0,1)$, $T_*$, and  $L_0>1$   such that
\eqref{3.14}-\eqref{3.19} hold simultaneously. Hence one can apply Remark \ref{remark2.2}(ii) to obtain the desired result. The proof of Lemma \ref{lem3.4} is now complete.
\end{proof}

By virtue of the above upper  and  lower solutions, we now prove Proposition \ref{prop3.1}.
\smallskip

\textbf{Proof of Proposition \ref{prop3.1}}.
For $t>T:=\max\{T_{\ast},T^{\ast}\}$, it follows from Lemmas  \ref{lem3.2}  and \ref{lem3.4}  that
$$\underline{h}(t)-s_{\mu,\rho}t\leq h(t)-s_{\mu,\rho}t\leq\overline{h}(t)-s_{\mu,\rho}t,$$
which leads to
$$-s_{\mu,\rho}T_{\ast}+L_0-1\leq h(t)-s_{\mu,\rho}t\leq -s_{\mu,\rho}T^{\ast}+\sigma+h(T^{\ast})+X_{0}.$$

Denote
$$C^{\ast}=\max\Big\{|-s_{\mu,\rho}T_{\ast}+L_0-1|,
|-s_{\mu,\rho}T^{\ast}+\sigma+h(T^{\ast})+X_{0}|,
\mathop{\max}_{t \in [0,T]}{|h(t)-s_{\mu,\rho}t|}\Big\}.$$
Then
$$|h(t)-s_{\mu,\rho}t|\le C^{\ast}~\textrm{for~all}~t\ge0.$$
Replacing the initial function of system \eqref{1.2} as $(u_{0}(-x),v_{0}(-x))$ and repeating the same arguments as above, one can find a positive constant $C_{\ast}$ such that $|g(t)+s_{\mu,\rho}t|\le  C_{\ast}$ for all $t\ge0$. Therefore, Proposition \ref{prop3.1} holds by taking $C=\max\{C^{\ast},C_{\ast}\}$. This completes the proof.
$\hfill\square$

\section{Sharp profile}\label{section 5}
\setcounter{equation}{0}
This section is devoted to proving Theorem  \ref{theo1.4}. As the proof is quite complicated, it will be divided into three subsections. In Subsection \ref{subsection 5.1}, we show that $h(t+t_{n})-s_{\mu,\rho}(t+t_{n})+2C\rightarrow L(t)$ and $(u,v)(x+s_{\mu,\rho}(t+t_{n})-2C,t+t_{n})\rightarrow(\hat{H},\hat{M})(x,t)$ as $n\to+\infty$ locally uniformly along a time sequence $\{t_{n}\}_{n=1}^{\infty}$ with $\lim\limits_{n\to+\infty}t_{n}=+\infty$. In Subsection \ref{subsection 5.2}, we demonstrate that $L(t)\equiv L(0)$ and $(\hat{H},\hat{M})(x,t)=(\phi_{s_{\mu,\rho}},\psi_{s_{\mu,\rho}})(L(0)-x)$. The proof of Theorem $\ref{theo1.4}$ is given in Subsection \ref{subsection 5.3}.

\subsection{Limit along a sequence $t_{n}\rightarrow+\infty$}\label{subsection 5.1}
By Proposition  \ref{prop3.1}, there exists $C>0$ such that
\begin{equation}\label{eq**1}
-C\le h(t)-s_{\mu,\rho}t\le C\text{ for }t\ge0.
\end{equation}

Define
$$H(x,t):=u(r(t)+x,t),\ \ M(x,t):=v(r(t)+x,t),\ \ l(t):=h(t)-r(t),$$
where
$$r(t):=s_{\mu,\rho}t-2C.$$
Obviously, $C\leq l(t)\leq3C$ for $t\ge0$ by \eqref{eq**1} and $(H,M,l)$ satisfies
\begin{equation*}
\begin{cases}
H_{t}=d_{1}H_{xx}+s_{\mu,\rho}H_{x}-aH+bM-F(H), &~g(t)-r(t)<x<l(t),~t>0,\\
M_{t}=d_{2}M_{xx}+s_{\mu,\rho}M_{x}+cH-dM-G(M), &~g(t)-r(t)<x<l(t),~t>0,\\
H(x,t)=M(x,t)=0, &~x=g(t)-r(t)~\textrm{or}~x=l(t),~t>0,\\
l'(t)=-\mu[M_{x}(l(t),t)+\rho H_{x}(l(t),t)]-s_{\mu,\rho}, &~t>0.
\end{cases}
\end{equation*}

Let $\{t_{n}\}_{n=1}^{\infty}$ be a positive sequence satisfying $\lim\limits_{n\to+\infty}t_{n}=+\infty$ and $\lim\limits_{n\to+\infty}l(t_{n})=\liminf\limits_{t\to+\infty}l(t)$. Define
$$(H_{n},M_{n})(x,t):=(H,M)(x,t+t_{n}),\ \ (l_{n},r_{n},g_{n})(t):=(l,r,g)(t+t_{n}).$$
\begin{lem}\label{lem4.1}
Subject to a subsequence,  we have
$$(H_{n},M_{n},l_n)\to (\hat{H},\hat{M},L)\text{ in }\Big[C^{1+\alpha,\frac{1+\alpha}{2}}_{\mathrm{loc}}(\Omega)\Big]^2\times C^{1+\frac{\alpha}{2}}_{\mathrm{loc}}(\mathbb{R})\text{ as }n\rightarrow+\infty,$$
where $\alpha \in (0,1)$ and  $\Omega:=\{(x,t): -\infty<x\leq L(t),~t \in \mathbb{R}\}$. Moreover, $(\hat{H},\hat{M},L)$ satisfies
\begin{equation}
\begin{cases}\label{4.1}
\hat{H}_{t}=d_{1}\hat{H}_{xx}+s_{\mu,\rho}\hat{H}_{x}-a\hat{H}+b\hat{M}-F(\hat{H}), &~(x,t) \in \Omega,\\
\hat{M}_{t}=d_{2}\hat{M}_{xx}+s_{\mu,\rho}\hat{M}_{x}+c\hat{H}-d\hat{M}-G(\hat{M}), &~(x,t) \in \Omega,\\
\hat{H}(L(t),t)=\hat{M}(L(t),t)=0, &~t \in \mathbb{R},\\
L'(t)=-\mu[\hat{M}_{x}(L(t),t)+\rho\hat{H}_{x}(L(t),t)]-s_{\mu,\rho},~L(t)\geq L(0), &~t \in \mathbb{R}.
\end{cases}
\end{equation}
\end{lem}
\begin{proof}
By Theorem  \ref{lem2.0}, there exists $C_{0}>0$ such that $0<h'(t)\leq C_{0}$ for $ t>0$, which yields
$$-s_{\mu,\rho}<l'_{n}(t)\leq C_{0}-s_{\mu,\rho}\text{ for }t>-t_{n}.$$
Denote
$$\xi=\frac{x}{l_{n}(t)},~(\tilde{H}_{n}(\xi,t),\tilde{M}_{n}(\xi,t))=(H_{n}(x,t),M_{n}(x,t)),
~\Sigma=\left(\frac{g_{n}(t)-r_{n}(t)}{l_{n}(t)},1\right).$$
Then $(\tilde{H}_{n},\tilde{M}_{n},l_{n})$~satisfies
\begin{equation}\label{eq4.2}
\begin{cases}
(\tilde{H}_{n})_{t}=\frac{d_{1}}{l_{n}^{2}(t)}(\tilde{H}_{n})_{\xi\xi}+\frac{\xi l_{n}'(t)+s_{\mu,\rho}}{l_{n}(t)}(\tilde{H}_{n})_{\xi}-a\tilde{H}_{n}+b\tilde{M}_{n}-F(\tilde{H}_{n}), &~\xi \in\Sigma,~t>-t_{n},\\
(\tilde{M}_{n})_{t}=\frac{d_{2}}{l_{n}^{2}(t)}(\tilde{M}_{n})_{\xi\xi}+\frac{\xi l_{n}'(t)+s_{\mu,\rho}}{l_{n}(t)}(\tilde{M}_{n})_{\xi}+c\tilde{H}_{n}-d\tilde{M}_{n}-G(\tilde{M}_{n}), &~\xi \in \Sigma,~t>-t_{n},\\
\tilde{H}_{n}(1,t)=\tilde{M}_{n}(1,t)=0, &~t>-t_{n},\\
l_{n}'(t)=-\frac{\mu}{l_{n}(t)}[(\tilde{M}_{n})_{\xi}(1,t)+\rho(\tilde{H}_{n})_{\xi}(1,t)]-s_{\mu,\rho}, &~t>-t_{n}.
\end{cases}
\end{equation}
According to to Theorem \ref{lem2.0}, $(u,v)$ is uniformly bounded for $x \in [g(t),h(t)]$ and $t>0$, which implies that $(H_{n},M_{n})$ is uniformly bounded for $x\le l_n(t)$ and $t\ge -t_n$. Fix $R>0$ and $T\in\mathbb{R}$. Then for any $p>1$, we can use the parabolic $L^p$ estimates cite{L1996} to \eqref{eq4.2} over $[-R-2,1]\times[T-2,T+1]$ to obtain that
$$\|(\tilde{H}_{n},\tilde{M}_{n})\|_{\big[W_{p}^{2,1}([-R-1,1]\times[T-1,T+1])\big]^2}\leq C_{R}~\textrm{for~all~large}~n,$$
where $C_{R}$ is a constant depending on $R$ and $p$ but independent of $n$ and $T$. Furthermore, for any $\alpha' \in (0,1)$, we can take $p>1$ large enough and use the Sobolev embedding theorem \cite{LSU1968} to deduce that
\begin{equation*}
\|(\tilde{H}_{n},\tilde{M}_{n})\|_{\big[C^{1+\alpha',\frac{1+\alpha'}{2}}([-R,1]\times
[T,+\infty))\big]^2}\leq \tilde{C}_{R}~\textrm{for~all~large}~n,
\end{equation*}
where $\tilde{C}_{R}$ is a constant depending on $R$ and $\alpha'$ but independent of $n$ and $T$. Combining with  \eqref{eq4.2}, we obtain
$$\|l_{n}\|_{C^{1+\frac{\alpha'}{2}}([T,+\infty))}\leq \hat{C}_{1}~\textrm{for~all~large}~n,$$
where $\hat{C}_{1}$ is a constant depending on $R$ and $\alpha'$ but independent of $n$ and $T$. Thus, by passing to a subsequence, still denoted by itself, we have
$$(\tilde{H}_{n},\tilde{M}_{n})\rightarrow (\overline{H},\overline{M})~\textrm{in}~ \Big[C^{1+\alpha,\frac{1+\alpha}{2}}_{\mathrm{loc}}((-\infty,1]\times\mathbb{R})\Big]^{2},~l_{n}\rightarrow L~\textrm{in}~ C^{1+\frac{\alpha}{2}}_{\mathrm{loc}}(\mathbb{R})$$
for some $\alpha \in (0,\alpha')$.
Applying the standard regularity theory to \eqref{eq4.2}, it follows that $(\overline{H},\overline{M},L)$ satisfies the following equations in the classical sense:
\begin{equation*}
\begin{cases}
\overline{H}_{t}=\frac{d_{1}}{L^{2}(t)}\overline{H}_{\xi\xi}+\frac{\xi L'(t)+s_{\mu,\rho}}{L(t)}\overline{H}_{\xi}-a\overline{H}+b\overline{M}-F(\overline{H}), &~\xi \in (-\infty,1],~t \in \mathbb{R},\\[2mm]
\overline{M}_{t}=\frac{d_{2}}{L^{2}(t)}\overline{M}_{\xi\xi}+\frac{\xi L'(t)+s_{\mu,\rho}}{L(t)}\overline{M}_{\xi}+c\overline{H}-d\overline{M}-G(\overline{M}), &~\xi \in (-\infty,1],~t \in \mathbb{R},\\[2mm]
\overline{H}(1,t)=\overline{M}(1,t)=0, &~t \in \mathbb{R},\\[2mm]
L'(t)=-\frac{\mu}{L(t)}\big[\overline{M}_{\xi}(1,t)+\rho\overline{H}_{\xi}(1,t)\big]-s_{\mu,\rho}, &~t \in \mathbb{R}.
\end{cases}
\end{equation*}
Set $(\hat{H}(x,t),\hat{M}(x,t))=(\overline{H}(x/{L(t)},t),\overline{M}(x/{L(t)},t))$. Then $(\hat{H},\hat{M},L)$ satisfies \eqref{4.1} and
$$\mathop{\lim}_{n\rightarrow+\infty}{\|{(H_{n},M_{n})-(\hat{H},\hat{M})}\|}
_{\big[C^{1+\alpha,\frac{1+\alpha}{2}}_{\mathrm{loc}}(\Omega)\big]^2}=~0.$$
Finally, since $L(t)=\lim\limits_{n\to+\infty}l(t+t_{n})$ and $L(0)=\lim\limits_{n\to+\infty}l(t_{n})=\liminf\limits_{t\to+\infty}l(t)$,  we obtain $L(t)\geq L(0)$ for $t\in \mathbb{R}$. This completes the proof.
\end{proof}
\subsection{Determining the limit pair $(\hat{H},\hat{M},L)$.} \label{subsection 5.2}
This subsection aims to demonstrate that
$$L(t)\equiv L(0)\text{ and }(\hat{H},\hat{M})(x,t)\equiv(\phi_{s_{\mu,\rho}},\psi_{s_{\mu,\rho}})(L(0)-x)\text{ for }(x,t) \in (-\infty,L(0)]\times\mathbb{R}.$$

In view of $C\leq l(t)\leq3C$ for $t\ge0$, we have
$$C\leq L(t)\leq3C\text{ for }t \in\mathbb{R}.$$
According to Lemma \ref{lem3.4},  for $x\in[-\underline{h}(t+t_{n})-r(t+t_{n}),\underline{h}(t+t_{n})-r(t+t_{n})]$ and $t\geq T_{\ast}-t_{n}$, there holds
\begin{equation}\label{4.4}
H_{n}(x,t)\geq(1-\tilde{\epsilon}e^{-\delta(t+t_{n}-T_{\ast})})\Phi_{n}(x,t),~M_{n}(x,t)\geq(1-\tilde{\epsilon}e^{-\delta(t+t_{n}-T_{\ast})})\Psi_{n}(x,t),
\end{equation}
 where $\underline{h}$ is defined in \eqref{lower sol} and
\begin{align*}
\begin{cases}
\Phi_{n}(x,t):=\phi_{s_{\mu,\rho}}(\underline{h}(t+t_{n})-r(t+t_{n})-x)+\phi_{s_{\mu,\rho}}(\underline{h}(t+t_{n})+r(t+t_{n})+x)-\phi_{s_{\mu,\rho}}(2\underline{h}(t+t_{n})),\\
\Psi_{n}(x,t):=\psi_{s_{\mu,\rho}}(\underline{h}(t+t_{n})-r(t+t_{n})-x)+\psi_{s_{\mu,\rho}}(\underline{h}(t+t_{n})+r(t+t_{n})+x)-\psi_{s_{\mu,\rho}}(2\underline{h}(t+t_{n})).
\end{cases}
\end{align*}
It follows from the definitions of $\underline{h}(t)$, $r(t)$, and $C$ that for $t+t_{n}\geq T_{\ast},$
$$\underline{h}(t+t_{n})-r(t+t_{n})
=2C-s_{\mu,\rho}T_{\ast}+L_0-(1-e^{-\delta(t+t_{n}-T_{\ast})})\ge2C-|-s_{\mu,\rho}T_{\ast}+L_0-1|\geq C $$
and
$$\underline{h}(t+t_{n})+r(t+t_{n})\rightarrow+\infty~\textrm{and}~\underline{h}(t+t_{n})\rightarrow+\infty~\textrm{as}~n\rightarrow+\infty,$$
which implies that for $x\leq C\leq L(t)$ and $t\in \mathbb{R}$,
\begin{equation*}
\begin{split}
&\liminf_{n\rightarrow+\infty}\Phi_{n}(x,t)
=\liminf_{n\rightarrow+\infty}\phi_{s_{\mu,\rho}}(\underline{h}(t+t_{n})-r(t+t_{n})-x)
\geq\phi_{s_{\mu,\rho}}(C-x),\\
&\liminf_{n\rightarrow+\infty}\Psi_{n}(x,t)
=\liminf_{n\rightarrow+\infty}\psi_{s_{\mu,\rho}}(\underline{h}(t+t_{n})-r(t+t_{n})-x)
\geq\psi_{s_{\mu,\rho}}(C-x).
\end{split}
\end{equation*}
Hence letting $n\rightarrow+\infty$ in \eqref{4.4}, we obtain
\begin{equation}\label{4.5}
(\hat{H},\hat{M})(x,t)\succeq(\phi_{s_{\mu,\rho}},\psi_{s_{\mu,\rho}})(C-x)\text{ for }x\leq C,~t\in\mathbb{R}.
\end{equation}

Define
\begin{equation*}
R^{\ast}:=\sup\left\{R\in\mathbb{R}:(\hat{H},\hat{M})(x,t)\succeq(\phi_{s_{\mu,\rho}},\psi_{s_{\mu,\rho}})(R-x)\text{ for }(x,t) \in (-\infty,R]\times\mathbb{R}\right\}.
\end{equation*}
By \eqref{4.5} and the fact that $(\hat{H},\hat{M})(L(t),t)=(0,0)$ with $C\le L(t)\le 3C$ for $t\in\mathbb{R}$, it follows that $R^{\ast}$ is well-defined,
\begin{equation}\label{4.6}
(\hat{H},\hat{M})(x,t)\succeq(\phi_{s_{\mu,\rho}},\psi_{s_{\mu,\rho}})(R^{\ast}-x)\text{ for }(x,t) \in(-\infty, R^{\ast}]\times\mathbb{R},
\end{equation}
and
\begin{equation}\label{4.7}
C\leq R^{\ast}\leq\mathop{\inf}_{t \in\mathbb{R}} L(t)=L(0).
\end{equation}
\begin{lem}\label{lem4.2}
$R^{\ast}=\mathop{\inf}_{t \in\mathbb{R}}L(t).$
\end{lem}
\begin{proof}
By \eqref{4.7}, it suffices to prove $R^{\ast}\geq\mathop{\inf}_{t \in\mathbb{R}} L(t)$. By contradiction, suppose  $R^{\ast}<\mathop{\inf}_{t \in\mathbb{R}}L(t)$. Then there exists a constant $\tau>0$ small enough such that
\begin{equation*}\label{Q1}
L(t)>R^{\ast}+\tau~\textrm{for~all}~t\in\mathbb{R}.
\end{equation*}
For clarity, the following arguments are divided into three steps.\\
\indent \textbf{Step 1.} We first show that
\begin{equation}\label{4.8}
(\hat{H},\hat{M})(x,t)\succ(\phi_{s_{\mu,\rho}},\psi_{s_{\mu,\rho}})(R^{\ast}-x)\text{ for }(x,t) \in(-\infty, R^{\ast}]\times\mathbb{R}.
\end{equation}

Indeed, if \eqref{4.8} fails, then \eqref{4.6} implies that there exists $(x_{0},t_{0}) \in (-\infty, R^{\ast}]\times\mathbb{R}$ such that
$$\hat{H}(x_{0},t_{0})=\phi_{s_{\mu,\rho}}(R^{\ast}-x_{0})~\textrm{or}~ \hat{M}(x_{0},t_{0})=\psi_{s_{\mu,\rho}}(R^{\ast}-x_{0}).$$
Without loss of generality, we assume that $\hat{H}(x_{0},t_{0})=\phi_{s_{\mu,\rho}}(R^{\ast}-x_{0})$. Recall that $(\hat{H},\hat{M})$ satisfies
\begin{equation}
\begin{cases}\label{4.9}
\hat{H}_{t}=d_{1}\hat{H}_{xx}+s_{\mu,\rho}\hat{H}_{x}-a\hat{H}+b\hat{M}-F(\hat{H}), &~(x,t) \in (-\infty,L(t)]\times\mathbb{R},\\
\hat{M}_{t}=d_{2}\hat{M}_{xx}+s_{\mu,\rho}\hat{M}_{x}+c\hat{H}-d\hat{M}-G(\hat{M}), &~(x,t) \in (-\infty,L(t)]\times\mathbb{R}.
\end{cases}
\end{equation}
It is easy to check that $(\phi_{s_{\mu,\rho}},\psi_{s_{\mu,\rho}})(R^{\ast}-x)$ also satisfies $(\ref{4.9})$ for $(x,t) \in(-\infty, R^{\ast}]\times\mathbb{R}$. Set
$$\vartheta(x,t)=\phi_{s_{\mu,\rho}}(R^{\ast}-x)-\hat{H}(x,t)\text{ for }(x,t) \in(-\infty, R^{\ast}]\times\mathbb{R}.$$
Obviously, $\vartheta(x_{0},t_{0})=0$ and \eqref{4.6} yields  $\vartheta(x,t)\leq0$ in $(-\infty, R^{\ast}]\times\mathbb{R}$.
From Theorem \ref{lem2.0}, there exists $C_0>0$ such that $\hat{H}(x,t)\le C_0$ in $(-\infty,R^*]\times\mathbb{R}$. Hence it follows from \eqref{4.6} and $\mathrm{(A1)}$ that there exists $\overline{K}>a+\max_{\xi\in[0,C_0]}F'(\xi)$ such that for $(x,t)\in(-\infty,R^*]\times\mathbb{R}$,
\begin{equation*}
\begin{split}
&\quad\ \vartheta_{t}-d_{1}\vartheta_{xx}-s_{\mu,\rho}\vartheta_{x}+\overline{K}\vartheta\\
&= a\hat{H}-b\hat{M}+F(\hat{H})-a\phi_{s_{\mu,\rho}}+b\psi_{s_{\mu,\rho}}-F(\phi_{s_{\mu,\rho}})+\overline{K}\vartheta\\
&\leq a\hat{H}+F(\hat{H})-a\phi_{s_{\mu,\rho}}-F(\phi_{s_{\mu,\rho}})+\overline{K}\vartheta\\
&\le-a\vartheta-\max_{\xi\in[0,C_0]}F'(\xi)\vartheta+\overline{K}\vartheta\le0.
\end{split}
\end{equation*}
Then by the strong maximum principle, we have  $\vartheta(x,t)\equiv0$ (i.e., $\phi_{s_{\mu,\rho}}(R^{\ast}-x)\equiv\hat{H}(x,t)$) for $(x,t)\in(-\infty, R^{\ast}]\times(-\infty,t_0]$, which yields  $\hat{H}(R^*,t_0)=\phi_{s_{\mu,\rho}}(0)=0$.
Applying the strong maximum principle to \eqref{4.1}, it is easy to deduce that
\begin{equation}\label{add equation 111}
\hat{H}(x,t)>0,\ \hat{M}(x,t)>0\text{ for }(x,t)\in(-\infty,L(t))\times\mathbb{R}.
\end{equation}
In particular,  $\hat{H}(R^*,t_0)>0$ due to $R^*<L(t)$ for $t\in\mathbb{R}$. We get a contradiction. Hence \eqref{4.8} holds.

\textbf{Step 2.} For $x\le R^*+\tau$, we define
\begin{equation*}
\begin{cases}
\omega_{1}(x):=\mathop{\sup}_{t \in \mathbb{R}}{[\phi_{s_{\mu,\rho}}(R^{\ast}-x)-\hat{H}(x,t)]},\\
\omega_{2}(x):=\mathop{\sup}_{t \in \mathbb{R}}{[\psi_{s_{\mu,\rho}}(R^{\ast}-x)-\hat{M}(x,t)]},
\end{cases}
\end{equation*}
where we understand that $(\phi_{s_{\mu,\rho}},\psi_{s_{\mu,\rho}})(R^{\ast}-x)$ is extended with $(\phi_{s_{\mu,\rho}},\psi_{s_{\mu,\rho}})(R^{\ast}-x)=(0,0)$ for $x\in(R^*,R^*+\tau]$.
We claim that
\begin{equation}\label{4.11}
\omega_{i}(x)<0\text{ for }x\le R^*+\tau \text{ and }i=1,2.
\end{equation}

Indeed,  $\omega_{i}(x)\leq0$ for $x\leq R^{\ast}$ and  $i=1,2$  by \eqref{4.6}. In light of $R^*+\tau<L(t)$ for $t\in\mathbb{R}$,   it follows from \eqref{add equation 111} that $\omega_{i}(x)\leq0$ for $x\in(R^*, R^{\ast}+\tau]$ and  $i=1,2$.
By contradiction, suppose  \eqref{4.11} fails. Then there exists $x_{0}\le R^{\ast}+\tau$ such that $\omega_{1}(x_{0})=0$ or $\omega_{2}(x_{0})=0$.
Without loss of generality, we assume that $\omega_{1}(x_{0})=0$. Clearly, it follows  from  \eqref{4.8} and \eqref{add equation 111} that  $\omega_{1}(x_{0})$ can not achieve zero at any finite time $t$. Thus, there exists a sequence $\{s_n\}_{n=1}^{\infty}\subset\mathbb{R}$ with $\lim\limits_{n\to+\infty}|s_{n}|=+\infty$ such that
\begin{equation}\label{Q2}
\mathop{\lim}_{n\rightarrow+\infty}\big[\hat{H}(x_0,s_{n})-\phi_{s_{\mu,\rho}}(R^{\ast}-x_0)\big]=0.
\end{equation}
Set
$$(\hat{H}_{n},\hat{M}_{n})(x,t)=(\hat{H},\hat{M})(x,t+s_{n}),~L_{n}(t)=L(t+s_{n}).$$
Repeating the same arguments used in the proof of Lemma $\ref{lem4.1}$ and passing to a subsequence if necessary, we can deduce that
$$(\hat{H}_{n},\hat{M}_{n},L_{n})\rightarrow(\hat{H}^{\ast},\hat{M}^{\ast},L^{\ast})~\textrm{in}~
\Big[C^{1+\alpha,\frac{1+\alpha}{2}}_{\mathrm{loc}}(\Omega^{\ast})\Big]^{2}\times C^{1+\frac{\alpha}{2}}_{\mathrm{loc}}(\mathbb{R})\text{ as }n\to+\infty,$$
where $\alpha \in (0,1)$ and $\Omega^{\ast}=\{(x,t):x\leq L^*(t),~t \in \mathbb{R}\}$. Moreover, $(\hat{H}^{\ast},\hat{M}^{\ast},L^{\ast})$ satisfies
\begin{equation}
\begin{cases}
\hat{H}^{\ast}_{t}=d_{1}\hat{H}^{\ast}_{xx}+s_{\mu,\rho}\hat{H}^{\ast}_{x}-a\hat{H}^{\ast}+b\hat{M}^{\ast}
-F(\hat{H}^{\ast}), &~x<L^{\ast}(t),~t \in \mathbb{R},\\
\hat{M}^{\ast}_{t}=d_{2}\hat{M}^{\ast}_{xx}+s_{\mu,\rho}\hat{M}^{\ast}_{x}+c\hat{H}^{\ast}-d\hat{M}^{\ast}
-G(\hat{M}^{\ast}), &~x<L^{\ast}(t),~t \in \mathbb{R},\\
\hat{H}^{\ast}(L^{\ast}(t),t)=\hat{M}^{\ast}(L^{\ast}(t),t)=0, &~t \in \mathbb{R},\\
(L^{\ast})'(t)=-\mu[\hat{M}^{\ast}_{x}(L^{\ast}(t),t)+\rho\hat{H}^{\ast}_{x}(L^{\ast}(t),t)]-s_{\mu,\rho}, &~t \in \mathbb{R}.
\end{cases}\label{4.12}
\end{equation}
 Moreover, \eqref{Q2} yields
 \begin{equation}\label{q eq *1}
 \hat{H}^*(x_0,0)=\phi_{s_{\mu,\rho}}(R^*-x_0)\text{ with }x_0\le R^*+\tau.
 \end{equation}

In the case when $x_0<R^*$, since $(\phi_{s_{\mu,\rho}},\psi_{s_{\mu,\rho}})(R^{\ast}-x)$ also satisfies \eqref{4.12} with $L^{\ast}(t)$ being replaced by $R^{\ast}$, similar to the arguments used in Step 1, we conclude from \eqref{q eq *1} and the strong maximum principle  that
$\hat{H}^{\ast}(x,t)\equiv\phi_{s_{\mu,\rho}}(R^{\ast}-x)$ for $(x,t)\in(-\infty,R^*]\times(-\infty,0]$.
In particular, we have $\hat{H}^*(R^*,0)=\phi_{s_{\mu,\rho}}(0)=0$, which is impossible since $L^*(0)>R^*$.

In the case when $x_0\in[R^*,R^*+\tau]$, it follows from  \eqref{q eq *1} that $\hat{H}^*(x_0,0)=0$. Similarly, we can apply the strong maximum principle to obtain that $\hat{H}^*(x,t)\equiv 0$ for $(x,t)\in(-\infty, L^*(t)]\times(-\infty,0]$, which yields $\hat{H}^*(R^*,0)=0$.  This is impossible since $L^*(0)>R^*$.  Therefore, \eqref{4.11} is valid.

\textbf{Step 3.} Completion of the proof.

In view of $\lim\limits_{x\to-\infty}(\phi_{s_{\mu,\rho}},\psi_{s_{\mu,\rho}})(R^{\ast}-x)=(u^*,v^*)$,  for any small $\epsilon_{0}>0$, there exists  $R_{0}=R_{0}(\epsilon_{0})<R^*$  such that
$$(\phi_{s_{\mu,\rho}},\psi_{s_{\mu,\rho}})(R^{\ast}-x)\succeq(u^*-\epsilon_{0},u^*-\epsilon_{0})\text{ for } x\leq R_{0}.$$
By \eqref{4.11} and the monotonicity of $\phi_{s_{\mu,\rho}}$ and  $\psi_{s_{\mu,\rho}}$, one can  choose $\epsilon \in (0,\epsilon_{0})$ small enough such that
\begin{equation}\label{4.13'}
(\phi_{s_{\mu,\rho}},\psi_{s_{\mu,\rho}})(R^{\ast}-R_{0}+\epsilon)\preceq(\phi_{s_{\mu,\rho}},\psi_{s_{\mu,\rho}})(R^{\ast}-R_{0})-(\omega_{1},\omega_{2})(R_{0}).
\end{equation}

Consider the following auxiliary problem
\begin{equation}
\begin{cases}
P_{t}=d_{1}P_{xx}+s_{\mu,\rho}P_{x}-aP+bQ-F(P), &~x<R_{0},~t>0,\\
Q_{t}=d_{2}Q_{xx}+s_{\mu,\rho}Q_{x}+cP-dQ-G(Q), &~x<R_{0},~t>0,\\
(P(R_{0},t),Q(R_{0},t))=(\phi_{s_{\mu,\rho}}(R^{\ast}-R_{0}+\epsilon),\psi_{s_{\mu,\rho}}(R^{\ast}-R_{0}+\epsilon)), &~t>0,\\
(P(x,0),Q(x,0))=(\phi_{s_{\mu,\rho}}(R^{\ast}-x),\psi_{s_{\mu,\rho}}(R^{\ast}-x)), &~x<R_{0}.
\end{cases}\label{4.14}
\end{equation}
Obviously, $(u^*,v^*)$ and $(\phi_{s_{\mu,\rho}},\psi_{s_{\mu,\rho}})(R^{\ast}-x)$ are a pair of upper and lower solutions of problem \eqref{4.14}. It follows from the comparison principle of cooperative systems that
\begin{equation*}\label{4.15}
(\phi_{s_{\mu,\rho}},\psi_{s_{\mu,\rho}})(R^{\ast}-x)\preceq(P,Q)(x,t)\preceq(u^*,v^*)\text{ for }x<R_0,~t>0,
\end{equation*}
 and $(P,Q)(x,t)$ is non-decreasing with respect to $t$ for each $x<R_0$. Hence there exist the limit functions
$$\lim\limits_{t\to+\infty}(P,Q)(x,t)=(P^{\ast},Q^{\ast})(x)\text{ for }x<R_{0},$$
where $(P^{\ast},Q^{\ast})$ satisfies
\begin{equation}
\begin{cases}
d_{1}P^{\ast}_{xx}+s_{\mu,\rho}P^{\ast}_{x}-aP^{\ast}+bQ^{\ast}-F(P^{\ast})=0, &~x<R_{0},\\
d_{2}Q^{\ast}_{xx}+s_{\mu,\rho}Q^{\ast}_{x}+cP^{\ast}-dQ^{\ast}-G(Q^{\ast})=0, &~x<R_{0},\\
(P^{\ast}(-\infty),Q^{\ast}(-\infty))=(u^*,v^*), \\
(P^{\ast}(R_{0}),Q^{\ast}(R_{0}))=(\phi_{s_{\mu,\rho}}(R^{\ast}-R_{0}+\epsilon),\psi_{s_{\mu,\rho}}(R^{\ast}-R_{0}+\epsilon)).
\end{cases}\label{4.16}
\end{equation}
Note that $(\phi_{s_{\mu,\rho}},\psi_{s_{\mu,\rho}})(R^{\ast}-x+\epsilon)$  satisfies \eqref{4.16} and it is an upper solution of problem \eqref{4.14}. Thus,
$$(P,Q)(x,t)\preceq(\phi_{s_{\mu,\rho}},\psi_{s_{\mu,\rho}})(R^{\ast}-x+\epsilon)\text{ for }x\leq R_{0},~t>0.$$
Letting $t\rightarrow+\infty$, we obtain
\begin{equation}\label{4.17'}
(P^{\ast},Q^{\ast})(x)\preceq(\phi_{s_{\mu,\rho}},\psi_{s_{\mu,\rho}})(R^{\ast}-x+\epsilon)\text{ for }x\leq R_{0}.
\end{equation}

We claim that
\begin{equation}\label{4.17}
(P^{\ast},Q^{\ast})(x)=(\phi_{s_{\mu,\rho}},\psi_{s_{\mu,\rho}})(R^{\ast}-x+\epsilon)\text{ for }x\leq R_{0}.
\end{equation}
Indeed, we define
$$\hat{P}(x):=P^{\ast}(x)-\phi_{s_{\mu,\rho}}(R^{\ast}-x+\epsilon)\text{ and }
~\hat{Q}(x):=Q^{\ast}(x)-\psi_{s_{\mu,\rho}}(R^{\ast}-x+\epsilon).$$
Then $(\hat{P},\hat{Q})$ satisfies
\begin{equation}
\begin{cases}
d_{1}\hat{P}_{xx}+s_{\mu,\rho}\hat{P}_{x}=a\hat{P}-b\hat{Q}+F(P^{\ast})-F(\phi_{s_{\mu,\rho}}), &~x<R_{0},\\
d_{2}\hat{Q}_{xx}+s_{\mu,\rho}\hat{Q}_{x}=-c\hat{P}+d\hat{Q}+G(Q^{\ast})-G(\psi_{s_{\mu,\rho}}), &~x<R_{0},\\
(\hat{P}(-\infty),\hat{Q}(-\infty))=(\hat{P}(R_{0}),\hat{Q}(R_{0}))=(0,0).
\end{cases}\label{4.18}
\end{equation}
Set
$$\hat{P}(\eta_{1}):=\mathop{\min}_{x \in (-\infty,R_{0}]}\hat{P}(x)~\text{ and }~
\hat{Q}(\eta_{2}):=\mathop{\min}_{x \in (-\infty,R_{0}]}\hat{Q}(x).$$
Then $\hat{P}(\eta_{1})\le0$ and $\hat{Q}(\eta_{2})\le0$ by \eqref{4.17'}. Thus, \eqref{4.17} is equivalent to $$\hat{P}(\eta_{1})=\hat{Q}(\eta_{2})=0.$$
By contradiction, suppose  $\hat{P}(\eta_{1})<0$. Then $P^{\ast}(\eta_{1})<\phi_{s_{\mu,\rho}}(R^{\ast}-\eta_{1}+\epsilon)$. Owing to $\hat{P}(R_0)=0$, we obtain $\eta_1<R_0$.  By the mean value theorem, there exists $\xi_{2} \in (P^{\ast}(\eta_{1}),\phi_{s_{\mu,\rho}}(R^{\ast}-\eta_{1}+\epsilon))$ such that
$$F(P^{\ast}(\eta_{1}))-F(\phi_{s_{\mu,\rho}}(R^{\ast}
-\eta_{1}+\epsilon))=F'(\xi_{2})\hat{P}(\eta_1).$$

If $a\hat{P}(\eta_{1})-b\hat{Q}(\eta_{2})+F'(\xi_{2})\hat{P}(\eta_{1})<0$, it then follows from  \eqref{4.18}  and the definition of $\hat{Q}(\eta_{2})$ that
$$0\leq d_{1}\hat{P}_{xx}(\eta_{1})+s_{\mu,\rho}\hat{P}_{x}(\eta_{1})\leq  a\hat{P}(\eta_{1})-b\hat{Q}(\eta_{2})+F'(\xi_{2})\hat{P}(\eta_{1})<0,$$
 a contradiction.

If $a\hat{P}(\eta_{1})-b\hat{Q}(\eta_{2})+F'(\xi_{2})\hat{P}(\eta_{1})\ge0$, then we have  $\hat{P}(\eta_{1})\ge b\hat{Q}(\eta_{2})/(a+F'(\xi_{2}))$ and  $\hat{Q}(\eta_{2})\le[a+F'(\xi_{2})]\hat{P}(\eta_{1})/b<0$, which yields $Q^*(\eta_2)<\psi_{s_{\mu,\rho}}(R^*-\eta_2+\epsilon)$. Thus, it follows from the mean value theorem that there exists $\xi_{3} \in (Q^{\ast}(\eta_{2}),\psi_{s_{\mu,\rho}}(R^{\ast}-\eta_{2}+\epsilon))$ such that
$$G(Q^{\ast}(\eta_{2}))-G(\psi_{s_{\mu,\rho}}(R^{\ast}-\eta_{2}+\epsilon))=G'(\xi_{3})\hat{Q}(\eta_2).$$
By \eqref{4.18}  and the definition of $\hat{P}(\eta_{1})$, we have
\begin{equation}\label{0}
\begin{split}
0&\leq d_{2}\hat{Q}_{xx}(\eta_{2})+s_{\mu,\rho}\hat{Q}_{x}(\eta_{2})\\
&\leq  -c\hat{P}(\eta_{1})+d\hat{Q}(\eta_{2})+G'(\xi_{3})\hat{Q}(\eta_{2})\\
&\leq \frac{\hat{Q}(\eta_{2})}{a+F'(\xi_{2})}[(a+F'(\xi_{2}))(d+G'(\xi_{3}))-bc],
\end{split}
\end{equation}
where the last inequality is due to $\hat{P}(\eta_{1})\ge b\hat{Q}(\eta_{2})/(a+F'(\xi_{2}))$.  Owing to $\mathrm{(A2)}$, there holds
$$\left(\frac{F(\zeta)}{\zeta}\right)'=\frac{F'(\zeta)\zeta-F(\zeta)}{\zeta^2}>0\text{ and }\left(\frac{G(\zeta)}{\zeta}\right)'=\frac{G'(\zeta)\zeta-G(\zeta)}{\zeta^2}>0\text{ for }\zeta>0.$$
In particular,   $F'(u^*)u^*>F(u^*)=-au^*+bv^*$ and $G'(v^*)v^*>G(v^*)=-dv^*+cu^*$,
which yields
\begin{equation}\label{**}
(a+F'(u^*))(d+G'(v^*))-bc>0.
\end{equation}

Recall that  $(P^{\ast}(-\infty),Q^{\ast}(-\infty))=(u^*,v^*)$  and
$\lim\limits_{x\to-\infty}(\phi_{s_{\mu,\rho}},\psi_{s_{\mu,\rho}})(R^{\ast}-x+\epsilon)=(u^*,v^*)$.  Noticing that  $\xi_{2} \in (P^{\ast}(\eta_{1}),\phi_{s_{\mu,\rho}}(R^{\ast}-\eta_{1}+\epsilon))$, $\xi_{3} \in (Q^{\ast}(\eta_{2}),\psi_{s_{\mu,\rho}}(R^{\ast}-\eta_{2}+\epsilon))$, and  $\eta_1,\eta_2\in(-\infty,R_0]$, we can choose $R_0<0$ smaller if necessary and use \eqref{**} to deduce that
$(a+F'(\xi_{2}))(d+G'(\xi_{3}))-bc>0$,
which leads to
$$\frac{\hat{Q}(\eta_{2})}{a+F'(\xi_{2})}[(a+F'(\xi_{2}))(d+G'(\xi_{3}))-bc]<0,$$
 contradicting \eqref{0}.

Therefore, it follows that $\hat{P}(\eta_{1})=0$  (i.e., $\hat{P}(x)=0$ for $x\leq R_{0}$).  Similar to the above arguments, we can prove that $\hat{Q}(x)=0$ for $x\leq R_{0}$. This proves \eqref{4.17}.

We are now ready to derive a contradiction by considering $(\hat{H},\hat{M})$, which  satisfies the first and second equations in \eqref{4.14}. In light of \eqref{4.8}  and \eqref{4.14}, we know that for any $t_* \in \mathbb{R}$ and $x<R_0$,
\begin{equation*}
(\hat{H},\hat{M})(x,t_*)\succ(\phi_{s_{\mu,\rho}},\psi_{s_{\mu,\rho}})(R^{\ast}-x)=(P,Q)(x,0).
\end{equation*}
Moreover, it  follows from  the definition of $\omega_i$, \eqref{4.13'},  and \eqref{4.14} that  for any $t_* \in \mathbb{R}$ and $t>0$,
\begin{equation*}
\begin{split}
&\hat{H}(R_{0},t+t_*)\geq\phi_{s_{\mu,\rho}}(R^{\ast}-R_{0})-\omega_{1}(R_{0})\ge
\phi_{s_{\mu,\rho}}(R^{\ast}-R_{0}+\epsilon)=P(R_{0},t),\\
&\hat{M}(R_{0},t+t_*)\geq\psi_{s_{\mu,\rho}}(R^{\ast}-R_{0})-\omega_{2}(R_{0})\ge
\psi_{s_{\mu,\rho}}(R^{\ast}-R_{0}+\epsilon)=Q(R_{0},t).
\end{split}
\end{equation*}
 Thus, $(\hat{H},\hat{M})(x,t+t_*)$ is an upper solution of problem \eqref{4.14} and
$$(\hat{H},\hat{M})(x,t+t_*)\succeq(P,Q)(x,t)~\textrm{for~all}~x\le R_{0},~t\ge 0,~t_* \in \mathbb{R},$$
which is equivalent to
$$(\hat{H},\hat{M})(x,t)\succeq(P,Q)(x,t-t_*)~\textrm{for~all}~x\le R_{0},~t\ge t_*,~t_* \in \mathbb{R}.$$
Letting $t_*\rightarrow-\infty$, we deduce from  \eqref{4.17} that
\begin{equation}\label{4.19}
(\hat{H},\hat{M})(x,t)\succeq(P^{\ast},Q^{\ast})(x)=(\phi_{s_{\mu,\rho}},
\psi_{s_{\mu,\rho}})(R^{\ast}-x+\epsilon)\text{ for }x\le R_{0},\ t \in \mathbb{R}.
\end{equation}

From \eqref{4.11} and the continuity of $\omega_i(x)$ in $x$,  there exists  $\epsilon_*>0$ such that $\omega_1(x), \omega_2(x)\le -\epsilon_*$ for $x\in[R_0,R^*+\tau]$. By the strict monotonicity of $\phi_{s_{\mu,\rho}}$ and $\psi_{s_{\mu,\rho}}$, we can choose $\epsilon_{1} \in (0,\min\{\tau,\epsilon_*\})$ small enough such that
\begin{equation}\label{4.19'}
(\phi_{s_{\mu,\rho}},\psi_{s_{\mu,\rho}})(R^{\ast}-x+\epsilon_{1})\preceq(\phi_{s_{\mu,\rho}},\psi_{s_{\mu,\rho}})(R^{\ast}-x)
+(\epsilon_*,\epsilon_*)\text{ for } x\in[R_{0},R^{\ast}+\epsilon_{1}],
\end{equation}
where we understand that $(\phi_{s_{\mu,\rho}},\psi_{s_{\mu,\rho}})(R^{\ast}-x)$ is extended with $(\phi_{s_{\mu,\rho}},\psi_{s_{\mu,\rho}})(R^{\ast}-x)=(0,0)$ for $x\in(R^*,R^*+\epsilon_1]$.
Hence it follows from \eqref{4.19'} that
for $x \in [R_{0},R^{\ast}+\epsilon_{1}]$ and $t \in \mathbb{R}$, there holds
\begin{align*}
&\quad(\hat{H},\hat{M})(x,t)-(\phi_{s_{\mu,\rho}},\psi_{s_{\mu,\rho}})(R^{\ast}-x+\epsilon_{1})\\
&\succeq (\hat{H},\hat{M})(x,t)-(\phi_{s_{\mu,\rho}},\psi_{s_{\mu,\rho}})(R^{\ast}-x)-(\epsilon_*,\epsilon_*)\\
&\succeq (\hat{H},\hat{M})(x,t)-(\phi_{s_{\mu,\rho}},\psi_{s_{\mu,\rho}})(R^{\ast}-x)
+(\omega_{1},\omega_{2})(x)\\
&\succeq(0,0),
\end{align*}
where the last inequality is due to the definitions of $\omega_1$ and $\omega_2$. Combining with $(\ref{4.19})$, we obtain
 $$(\hat{H},\hat{M})(x,t)\succeq(\phi_{s_{\mu,\rho}},\psi_{s_{\mu,\rho}})(R^{\ast}-x+\epsilon_{1})~\textrm{for}
 ~x\leq R^{\ast}+\epsilon_{1},~t \in \mathbb{R}.$$
This contradicts  the definition of $R^{\ast}$. The proof is complete.
\end{proof}
\begin{prop}\label{prop4.3}
$(\hat{H},\hat{M})(x,t)\equiv(\phi_{s_{\mu,\rho}},\psi_{s_{\mu,\rho}})(L(0)-x)$ and $L(t)\equiv L(0)$ for $x\le L(0)$ and $t\in\mathbb{R}$.
\end{prop}
\begin{proof}
It follows from Lemma \ref{lem4.2}  that $R^{\ast}=\inf_{t\in R}L(t)=L(0)$. Thus, there holds $$(\hat{H},\hat{M})(L(0),0)=(0,0)=(\phi_{s_{\mu,\rho}},\psi_{s_{\mu,\rho}})(R^{\ast}-L(0)).$$
Moreover, \eqref{4.6} yields that $(\hat{H},\hat{M})(x,t)\succeq(\phi_{s_{\mu,\rho}},\psi_{s_{\mu,\rho}})(R^{\ast}-x)$ for $(x,t)\in(-\infty,L(0)]\times\mathbb{R}$.
By contradiction, suppose $(\hat{H},\hat{M})(x,t)\not\equiv(\phi_{s_{\mu,\rho}},\psi_{s_{\mu,\rho}})(R^{\ast}-x)$ for $(x,t)\in(-\infty,L(0)]\times\mathbb{R}$. It then  follows from the strong maximum principle and the Hopf boundary lemma that
$$(\hat{H}_{x},\hat{M}_{x})(L(0),0)\prec-(\phi_{s_{\mu,\rho}}',\psi_{s_{\mu,\rho}}')(0),$$
which leads to
\begin{equation}\label{add 01}
\hat{M}_x(L(0),0)+\psi_{s_{\mu,\rho}}'(0)+\rho[\hat{H}_x(L(0),0)+\phi_{s_{\mu,\rho}}'(0)]<0.
\end{equation}
Since $L(0)=\inf_{t\in\mathbb{R}}L(t)$, we have $L'(0)=0$. Then \eqref{s mu} and the last equation in \eqref{4.1} imply that
$$\hat{M}_x(L(0),0)+\psi_{s_{\mu,\rho}}'(0)+\rho[\hat{H}_x(L(0),0)+\phi_{s_{\mu,\rho}}'(0)]=0,$$
 contradicting \eqref{add 01}. Hence, $(\hat{H},\hat{M})(x,t)\equiv(\phi_{s_{\mu,\rho}},\psi_{s_{\mu,\rho}})(R^{\ast}-x)$ for $(x,t)\in(-\infty,L(0)]\times\mathbb{R}$. This implies that $\phi_{s_{\mu,\rho}}(R^*-L(t))=\hat{H}(L(t),t)=0$ for $t\in\mathbb{R}$. Thus, it follows from the strict monotonicity of $\phi_{s_{\mu,\rho}}$ that $L(t)\equiv R^*=L(0)$ for $t\in\mathbb{R}$. This completes the proof.
\end{proof}
\subsection{Proof of Theorem  \ref{theo1.4}}\label{subsection 5.3}
We will complete the proof of Theorem  \ref{theo1.4}  in this subsection. For clarity, we  begin by proving two claims.

\textbf{Claim 1:} Let $\{t_n\}_{n=1}^{\infty}$ be the sequence in Lemma  \ref{lem4.1}. Then for any $t\in\mathbb{R}$, we have $\lim\limits_{n\to\infty}h'(t+t_n)=s_{\mu,\rho}$. Moreover, \eqref{1.9} and \eqref{1.10} hold along $t=t_n$.

It follows from Lemma \ref{lem4.1} and Proposition \ref{prop4.3} that
 $$\lim\limits_{n\to+\infty}l(t+t_{n})=\lim\limits_{n\to+\infty}\left[h(t+t_{n})-r(t+t_{n})\right]=L(0)\text{ in }C^{1+\frac{\alpha}{2}}_{\mathrm{loc}}(\mathbb{R}).$$
In view of $r(t)=s_{\mu,\rho}t-2C$, we obtain $\lim\limits_{n\to+\infty}h'(t+t_{n})=s_{\mu,\rho}$ in $C^{\frac{\alpha}{2}}_{\mathrm{loc}}(\mathbb{R})$. Thus, according to Lemma \ref{lem4.1}  and Proposition \ref{prop4.3}, one then finds that
$$\lim\limits_{n\to+\infty}(u,v)(x+h(t+t_{n}),t+t_{n})=(\phi_{s_{\mu,\rho}},\psi_{s_{\mu,\rho}})(-x)~\textrm{in}~
\Big[C^{1+\alpha,\frac{1+\alpha}{2}}_{\mathrm{loc}}((-\infty,0]\times\mathbb{R})\Big]^2.$$
In particular, for any $M_{0}>0$, there holds
\begin{equation}\label{11}
\mathop{\lim}_{n\rightarrow+\infty}\|(u,v)(x,t_{n})-(\phi_{s_{\mu,\rho}},\psi_{s_{\mu,\rho}})(h(t_{n})-x)\|
_{[L^{\infty}([h(t_{n})-M_{0},h(t_{n})])]^2}=0.
\end{equation}

On the other hand, for any given $\epsilon>0$ small, Lemmas \ref{lem3.2} and \ref{lem3.4} imply that there exist constants $M_{1}>0$ and   $N\geq1$ such that for $n\ge N$,
$$(u^*-\epsilon,v^*-\epsilon)\preceq(u,v)(x,t_{n})\preceq(u^*+\epsilon,v^*+\epsilon)~\mathrm{for}
~x\in[0,h(t_{n})-M_{1}].$$
Due to  $(\phi_{s_{\mu,\rho}}(+\infty),\psi_{s_{\mu,\rho}}(+\infty))=(u^*,v^*)$, there exists   $M_{2}>0$ large enough such that for $n\ge N$,
$$(u^*-\epsilon,v^*-\epsilon)\preceq(\phi_{s_{\mu,\rho}},\psi_{s_{\mu,\rho}})(h(t_{n})-x)
\preceq(u^*,v^*)~\mathrm{for}~x\in(-\infty,h(t_{n})-M_{2}].$$
Thus, taking $M_{0}=\max\{M_{1},M_{2}\}$, it follows that for $n\geq N$,
$$\|(u,v)(x,t_{n})-(\phi_{s_{\mu,\rho}},\psi_{s_{\mu,\rho}})(h(t_{n})-x)\|
_{[L^{\infty}([0,h(t_{n})-M_{0}])]^2}<2\epsilon.$$
By using the arbitrariness of $\epsilon$ and combining with \eqref{11}, we obtain
\begin{equation}\label{4.20}
\mathop{\lim}_{n\rightarrow+\infty}\|(u,v)(x,t_{n})-
(\phi_{s_{\mu,\rho}},\psi_{s_{\mu,\rho}})(h(t_{n})-x)\|
_{[L^{\infty}([0,h(t_{n})])]^2}=0.
\end{equation}
Similarly, one can consider system \eqref{1.2} with the initial function $(u_{0}(-x),v_{0}(-x))$ to deduce that
\begin{equation}\label{4.21}
\mathop{\lim}_{n\rightarrow+\infty}\|(u,v)(x,t_{n})-
(\phi_{s_{\mu,\rho}},\psi_{s_{\mu,\rho}})(h(t_{n})-x)\|
_{[L^{\infty}([g(t_{n}),0])]^2}=0.
\end{equation}

\textbf{Claim 2:} \emph{$\lim\limits_{t\to+\infty}(h(t)-s_{\mu,\rho}t)=h^{\ast}$, where $h^*:=L(0)-2C$.}

According to Claim 1, along a sequence $\{t_{n}\}$ satisfying
$$\lim\limits_{n\to+\infty}\big(h(t_{n})-s_{\mu,\rho}t_{n}+2C\big)=\liminf\limits_{t\to+\infty}\big(h(t)-s_{\mu,\rho}t+2C\big)
=L(0),$$
we know that \eqref{4.20} holds, and
\begin{equation}\label{4.22}
\lim\limits_{n\rightarrow+\infty}\big(h(t_{n})-s_{\mu,\rho}t_{n}\big)=\liminf\limits_{t\rightarrow+\infty}\big(h(t)-s_{\mu,\rho}t\big)=h^*,\ \ \lim\limits_{n\rightarrow+\infty}h'(t_{n})=s_{\mu,\rho}.
\end{equation}
By contradiction, suppose Claim 2 fails. Then there exists a constant $\tilde{h}^{\ast}>h^{\ast}$ and  a sequence $\{s_{n}\}_{n=1}^{\infty}$ with $\lim\limits_{n\to+\infty}s_n=+\infty$  such that
$$\lim\limits_{n\to+\infty}\big(h(s_{n})-s_{\mu,\rho}s_{n}\big)=\limsup\limits_{t\to+\infty}\big(h(t)-s_{\mu,\rho}t\big)=\tilde{h}^{\ast}.$$

We next demonstrate that $(\overline{u},\overline{v},\overline{g},\overline{h})$, which is defined in \eqref{3.1} with $X_{0}=\sigma=(\tilde{h}^{\ast}-h^{\ast})/4>0$ and $T^{\ast}=t_{n}$, is an upper solution to system \eqref{1.2}.  Clearly, for such $\sigma$, there exist  constants $\delta>0$ and $K_{1}>0$ such that \eqref{3.8}, \eqref{3.10}, and \eqref{3.11} hold simultaneously. Hence from the proof of Lemma \ref{lem3.2}, we know that  \eqref{3.2},  \eqref{3.3},  \eqref{3.4},  \eqref{3.5}, and  \eqref{3.6} hold.
It remains to prove  \eqref{3.7}.
Indeed, it follows from \eqref{4.20} and \eqref{4.21} that for
 large $n$ and $x\in[\overline{g}(t_{n}),h(t_{n})]=[g(t_{n}),h(t_{n})]$, there holds
 \begin{equation*}
\begin{split}
&\overline{u}(x,t_{n})=(1+K_{1})\phi_{s_{\mu,\rho}}(\overline{h}(t_{n})-x)=(1+K_{1})\phi_{s_{\mu,\rho}}(h(t_n)+X_0-x)
\geq u(x,t_{n}),\\
&\overline{v}(x,t_{n})=(1+K_{1})\psi_{s_{\mu,\rho}}(\overline{h}(t_{n})-x)=(1+K_{1})\psi_{s_{\mu,\rho}}(h(t_n)+X_0-x)
\geq v(x,t_{n}),
\end{split}
\end{equation*}
that is, \eqref{3.7} holds for large $n$. Therefore, it follows from Lemma \ref{lem2.1} that
\begin{equation*}\label{4.23}
h(t)\leq\overline{h}(t)~\mathrm{for}~t>t_{n}.
\end{equation*}
In particular, for all large integer $k$ satisfying $s_{k}\geq t_{n}$, we have
$$h(s_{k})\leq\overline{h}(s_{k})=s_{\mu,\rho}(s_{k}-t_{n})+\sigma(1-e^{-\delta(s_{k}-t_{n})})+h(t_{n})+X_{0},$$
which implies that
$$\tilde{h}^{\ast}=\lim_{k\rightarrow+\infty}(h(s_{k})-s_{\mu,\rho}s_{k})\leq-s_{\mu,\rho}t_{n}+\sigma+h(t_{n})+X_{0}.$$
Letting $n\rightarrow+\infty$, we obtain
$$\tilde{h}^{\ast}\leq h^{\ast}+\sigma+X_{0}=h^{\ast}+(\tilde{h}^{\ast}-h^{\ast})/2<\tilde{h}^{\ast},$$
a contradiction. Thus, Claim 2 holds

According to Claim 2, any positive sequence $\{t_n\}_{n=1}^{\infty}$ converging to $+\infty$ can be used in Lemma \ref{lem4.1}, and hence it has a sequence such that \eqref{4.20} and \eqref{4.22} hold. This clearly implies that \eqref{1.8} and \eqref{1.10} are valid.
Similarly, by considering  system \eqref{1.2} with the initial function $(u_{0}(-x),v_{0}(-x))$, one can prove that  \eqref{1.7} and \eqref{1.9} also hold. This accomplishes the proof of Theorem \ref{theo1.4}.
$\hfill\square$

\section{Discussions}\label{section 6}
\setcounter{equation}{0}
In \cite{MM1996},  Molina-Meyer proposed a reaction-diffusion cooperative system in a bounded domain, denoted by \eqref{Cauchy Problem}, and primarily investigated the properties of its positive steady states. From a biological perspective, such a system is valuable for modeling the invasion dynamics of two cooperative species. In \cite{YD2010}, Du and Lin extended the Fisher-KPP equation by introducing free boundaries to study  the invasion behaviors of a single species. Their work established key concepts such as the spreading-vanishing dichotomy and the asymptotic spreading speed of the free boundary. Compared to  fixed boundary problems, free boundary problems are more realistic as they can accurately locate the spreading fronts \cite{Du2022BMS}. Inspired by these seminal works \cite{Du2022BMS}, \cite{YD2010}, \cite{MM1996},   we introduce temporal changes in the living regions in cooperative system  \eqref{Cauchy Problem} and study the dynamics of the resulting system, denoted by \eqref{1.2}.

To begin with, we investigate the long-term behavior of the solution to system  \eqref{1.2}, with the results presented as a spreading-vanishing dichotomy in Theorem \ref{theo1.1}. In the vanishing case, the living region of the species is finite, and global convergence of the population densities is established (see Theorem \ref{theo1.1}(i)). In the spreading scenario, the living region of the species expands infinitely and the population densities converge locally uniformly to the unique positive equilibrium $(u^*,v^*)$   in $\mathbb{R}$ (see Theorem \ref{theo1.1}(ii)).

To gain a deeper understanding of species   invasion, we perform a rigorous analysis of the asymptotic spreading speeds of the moving fronts and the global convergence of  population densities over the habitat interval $(g(t),h(t))$.
The key findings are summarized in Theorems \ref{theo1.3} and \ref{theo1.4}. Firstly, as shown in  Theorem \ref{theo1.3}, in spreading cases, the two moving fronts $h(t)$ and $g(t)$ possess the same asymptotic spreading speed $s_{\mu,\rho}$, which is uniquely determined by the species'  invasion ability (see \eqref{s mu}). Specifically, the two moving fronts have  the following properties:
$h(t)=s_{\mu,\rho} t+o(t)$ and $g(t)=-s_{\mu,\rho} t+o(t)$ as $t\to+\infty$. Furthermore, Theorem \ref{theo1.4} (see   \eqref{1.7} and \eqref{1.8})  refines this understanding by specifying that
 $\lim\limits_{t\to+\infty}\big(h(t)-s_{\mu,\rho}t-h^*\big)
 =\lim\limits_{t\to+\infty}\big(g(t)+s_{\mu,\rho}t-g^*\big)=0$
for constants $h^*$ and $g^*$.
This indicates that, in the spreading scenario,  the moving fronts $h(t)$ and $g(t)$ ultimately converge to linear functions with constant correction terms.
Secondly,  the global behavior of the population densities $u(x,t)$ and $v(x,t)$ across the entire living area $(g(t),h(t))$ is established in Theorem \ref{theo1.4} (see \eqref{1.9} and \eqref{1.10}).
These results demonstrate that the long-term dynamics of the moving fronts and population densities can be characterized by the semi-wave solution derived from system \eqref{1.5}. It is  noteworthy that  the rigorous study of sharp asymptotic spreading speeds for fronts and the global profile for population densities in multi-species systems has received limited attention  in the literature.
Hence, the precise dynamical results provided in Theorem \ref{theo1.4} highlight their  significance in advancing our understanding of these complex ecological systems.

In conclusion, this paper investigates the  dynamics of a reaction-diffusion system with two free boundaries, modeling the invasion of two cooperative species. The main results cover the spreading-vanishing dichotomy, the asymptotic spreading speeds of the free boundaries, and the  sharp global profiles of the population densities  (see Theorems \ref{theo1.1}, \ref{theo1.3}, and \ref{theo1.4}). The approaches in this paper are based on the method of upper and lower solutions, the maximum principle, and a blowup argument. It is hoped that the analysis in this paper will inspire further research on cooperative systems with free boundaries.

\bigskip
\noindent\textbf{Acknowledgments:}
The authors are very grateful to the handling editor and anonymous referee for careful reading and valuable comments, which led to improvements of the manuscript.

\noindent\textbf{Declaration of competing interest:} Authors state no conflict of interest.

\noindent\textbf{Research funding:}  H. Nie was partially supported by the National Natural Science Foundation of China (12471464), the Natural Science Basic Research Program of Shaanxi (2023-JC-JQ-03), and the Fundamental Research Funds for the Central Universities (GK202402004).
Z. Wang was partially supported by the National Natural Science Foundation of China (12371496, 12471168).

\noindent\textbf{Authorship contribution statement:}
\textbf{Qian Qin}: Formal analysis, Writing-original draft $\&$ editing.
\textbf{Jinjing Jiao}: Formal analysis $\&$ Writing-original draft.
\textbf{Zhiguo Wang}: Conceptualization, Methodology, Supervision $\&$ Writing-review.
\textbf{Hua Nie}: Methodology, Supervision $\&$ Writing-review.

\noindent\textbf{Data availability:} No data was used for the research described in the article.


\end{document}